# APPROXIMATING A SEQUENCE OF OBSERVATIONS BY A SIMPLE PROCESS[1]

By Dinah Rosenberg, Eilon Solan[2] and Nicolas Vieille

*Université Paris Nord, Northwestern University and Tel Aviv University, and HEC*

Given an arbitrary long but finite sequence of observations from a finite set, we construct a simple process that approximates the sequence, in the sense that with high probability the empirical frequency, as well as the empirical one-step transitions along a realization from the approximating process, are close to that of the given sequence.

We generalize the result to the case where the one-step transitions are required to be in given polyhedra.

**1. Introduction.** In a seminal work, Baum and Petrie (1966) studied the following problem. Can one recover a homogenous hidden Markov chain from a finite sample $(x_0, x_1, \ldots, x_N)$ from the chain. They prove that the maximum likelihood estimate converges to the correct value, as $N$ goes to infinity.

This problem has several applications, including ecology [Baum and Eagon (1967)], speech recognition [see, e.g., Rabiner (1989)] and identifying gene structure [see, e.g., Krogh, Mian and Haussler (1994)].

We study the following related problem. Can one find a "simple" process $(s_n)$ that "explains" a given observation $(x_0, x_1, \ldots, x_N)$? More specifically, we are given a finite sequence $(x_0, \ldots, x_N)$ out of a finite set $S$, and we would like to find a simple $S$-valued process $(s_n)$ that satisfies the following two properties:

(i) under $(s_n)_{n \leq N}$, with high probability, the empirical frequency of $s \in S$ is close to the frequency of stages $m < N$ such that $x_m = s$, and

Received May 2002; revised November 2003.
[1]Supported by the Arc-en-Ciel/Keshet program for 2001/2002.
[2]Supported by Israel Science Foundation Grant No. 69/01-1.

(ii) the conditional law of $s_{n+1}$, given $(s_0, \ldots, s_n)$, is close to the empirical frequency of one-step transitions from $s_n$ to $s_{n+1}$ in $(x_0, \ldots, x_N)$ (i.e., the frequency of stages $m < N$ such that $x_{m+1} = s_{n+1}$ out of the stages $m < N$ such that $x_m = s_n$).

Were only the property (i) required, an i.i.d. sequence would do. The simplest processes that allow for serial correlation are homogeneous Markov chains. Therefore, a naive solution to this problem is to define $(s_n)$ to be the homogeneous Markov chain in which the transition from $s$ to $s'$ is the frequency of stages $m < N$ such that $x_{m+1} = s'$ out of the stages $m < N$ such that $x_m = s$. It is true that *asymptotically* this Markov chain satisfies our requirements. However, we wish to have an approximation at time $N$, where $N$ is the number of observations, and the naive Markov chain may fail to do so. The concept of a simple process we use is therefore slightly more complicated: a simple process in our context is a *piecewise* homogeneous Markov chain with a *bounded* number of pieces. Our basic result states that, provided $N$ is large enough, *every* sequence can be explained, in the above sense, by a piecewise homogeneous Markov chain with at most $|S|$ pieces. Our proof is constructive, in the sense that we provide an algorithm that produces the desired piecewise homogeneous Markov chain.

We also analyze a more general question. It is sometimes the case that the process we construct has to satisfy some exogenous constraints, for example, the one-step transitions must belong to some pre-defined polyhedra of probability measures. These polyhedra may reflect some a priori knowledge of the physics of the problem at hand. We then have to construct a process such that both (i) and (ii) are satisfied, and, in addition, the conditional law of $s_{n+1}$, given $(s_0, \ldots, s_n)$, must belong to some polyhedron $V(s_n)$. So that (ii) will hold, the empirical transitions along the observed sequence must be close to the given polyhedra. We prove that under proper conditions, and if $N$ is large enough, there exists a piecewise *hidden* Markov chain with a bounded number of pieces that satisfies these three requirements.

A consequence of our result is the following. Let $(\mathbf{z}_n)$ be any $S$-valued process such that the conditional law of $\mathbf{z}_{n+1}$, given $(\mathbf{z}_0, \ldots, \mathbf{z}_n)$, belongs to some given polyhedron $V(\mathbf{z}_n)$, a.s. for each $n < N$. Assume, moreover, that there is an irreducible transition function $b$ such that $b(s, \cdot) \in V(s)$ for every $s$. Then, provided $N$ is large enough, for most realizations $(x_0, x_1, \ldots, x_N)$ of $(\mathbf{z}_n)$ one can find a piecewise homogeneous hidden Markov chain $(\mathbf{s}_n)$ with at most $|S|$ pieces such that both (i) and (ii) above hold, and the conditional law of $\mathbf{s}_{n+1}$, given $(\mathbf{s}_0, \ldots, \mathbf{s}_n)$, belongs to $V(\mathbf{s}_n)$, a.s. for every $n < N$. More precisely, the measure of the set of realizations that can be explained in the sense we just described goes to 1 as $N$ goes to infinity. Thus, most realizations from $(\mathbf{z}_n)$ can be explained by a simple process. In other words, assuming *only* that the sequence $(x_0, \ldots, x_N)$ is generated by a process that satisfies the given physical constraints, it is very likely that a simple process can be found, which with high probability has the same empirical behavior as the given sequence $(x_0, \ldots, x_N)$.

The paper is organized as follows. In Section 2 we define and investigate the problem with no polyhedral restriction, and in Section 3 we turn to the general problem.

**2. The basic problem.** Given a finite set $K$, we let $|K|$ denote the number of elements in $K$, and $\mathcal{P}(K)$ denote the space of probability distributions over $K$. Throughout the paper we fix a finite set $S$ of *states*. We use the symbol "$\subset$" to denote strict inclusion. For every subset $C \subseteq S$, $\overline{C} = S \setminus C$ is the complement of $C$ in $S$.

2.1. *Presentation.* The basic problem can be stated as follows. A sequence $x = (x_0, x_1, \ldots, x_N)$ in $S$ with finite length $N + 1 \in \mathbf{N}$ is given. An *observer* gets to see this sequence, or at least gets to know the number $N^x_{s,t} = |\{n < N | (x_n, x_{n+1}) = (s,t)\}|$ of one-step transitions from $s$ to $t$, for each $s, t \in S$. The observer wishes to find a *simple* stochastic process $(\mathbf{z}_n)_n$ over $S$, such that any typical realization of $\mathbf{z}_0, \ldots, \mathbf{z}_N$ fits the data. We proceed to give a formal meaning to this question before we state our basic result.

For $s \in S$, denote by

$$N^x_s = \sum_{t \in S} N^x_{s,t} = |\{n < N : x_n = s\}| \quad \text{and} \quad \nu^x(s) = \frac{N^x_s}{N}$$

the *number of stages spent in $s$* along $x$ (excluding $x_N$), and the *observed occupancy measure* of $s$, respectively. The *observed transition function* $p^x$ is defined by

(1) $$p^x(s,t) = \frac{N^x_{s,t}}{N^x_s} \quad \text{for each } s, t \in S \text{ s.t. } N^x_s > 0.$$

If $N^x_s = 0$, the definition of $p^x(s, \cdot) \in \mathcal{P}(S)$ is irrelevant.

The most natural notion of a simple process is that of a homogeneous Markov chain. As is argued in Remark 1 below, this notion is not flexible enough to allow for a good approximation in finite time. Thus, we introduce the notion of a *piecewise homogeneous Markov chain*.

DEFINITION 1. Let $K$ be a positive integer. A process $\mathbf{z} = (\mathbf{z}_n)_{0 \leq n \leq N}$ is a *piecewise homogeneous Markov chain with $K$ pieces* if (i) $\mathbf{z}$ is a Markov chain and (ii) there exist integers $0 = n_0 \leq n_1 \leq \cdots \leq n_K = N$ such that the random variables $(\mathbf{z}_n), n_{k-1} \leq n < n_k$, form a homogeneous Markov chain, for each $k = 1, \ldots, K$.

The law of an $S$-valued process $\mathbf{z} = (\mathbf{z}_n)_{n \leq N}$ is denoted by $\mathbf{P_z}$. If $\mathbf{z}$ is a (possibly nonhomogeneous) Markov chain, we denote by $p_n^{\mathbf{z}}$, $n < N$, the conditional distribution of $\mathbf{z}_{n+1}$ given $\mathbf{z}_n$. Also, $\nu_m^{\mathbf{z}}$ is the *empirical occupancy measure* in the first $m$ stages:

$$\nu_m^{\mathbf{z}}(s) = \frac{1}{m}|\{0 \leq n \leq m-1 : \mathbf{z}_n = s\}|.$$

We are now in a position to state our basic theorem.

THEOREM 1. *For every $\varepsilon > 0$, every $\delta \in (0, \frac{1}{2(4|S|+1)})$ and every $\zeta \in (0, 2\delta)$, there exists $N_0 \in \mathbf{N}$ such that the following holds. For every $N \geq N_0$ and every $S$-valued sequence $x = (x_0, \ldots, x_N)$, there is a piecewise homogeneous Markov chain $\mathbf{z}$ over $S$, with at most $|S|$ pieces, such that:*

(B1) $\mathbf{P_z}(|\frac{\nu_N^{\mathbf{z}}(s)}{\nu^x(s)} - 1| \geq \varepsilon) \leq \frac{1}{N^\zeta}$ *for every $s \in S$ that satisfies $\nu^x(s) \geq \frac{1}{N^\delta}$.*
(B2) $\mathbf{P_z}$-*a.s., one has $\|p_n^{\mathbf{z}}(\mathbf{z}_n, \cdot) - p^x(\mathbf{z}_n, \cdot)\|_\infty \leq \varepsilon$ for each $n \neq n_k$.*

Leaving aside the technical qualifications, Theorem 1 has the following implications. The number of pieces of the approximating process is independent of the length of the sequence $x$. In all stages, with the possible exception of at most $|S|$ of them, the transition function of $\mathbf{z}$ is very close to the observed transition function $p^x$. Moreover, for any typical realization of the first $N$ components of $\mathbf{z}$, the empirical occupancy measure $\nu_N^{\mathbf{z}}$ is very close to the observed occupancy measure $\nu^x$ [restricted to states $s$ whose observed occupancy measure $\nu^x(s)$ is not negligible].

We stress that we consider realizations $\mathbf{z}$ of the *same* length as the sequence. In that sense, our result is not an asymptotic result, but provides the basis for a good approximation in finite time, provided the sequence is long enough. Our proof is constructive, in the sense that we provide an algorithm that can be used to construct $\mathbf{z}$.

Observe that (B1) and (B2) are not exactly of the same nature. Indeed, (B1) relates to the samples from $\mathbf{z}$, while (B2) is a structural property of $\mathbf{z}$. From the proof it will be clear that endless variations are possible.

REMARK 1. The naive solution is to consider a Markov chain $\mathbf{z}$ with transition function $p^x$. However, such a process may fail to yield a good approximation in finite time. Indeed, let $S = \{a, b\}$, and consider the sequence $x = (a, a, \ldots, a, b, b, \ldots, b, a)$ that contains $N$ $a$'s followed by $N$ $b$'s, and ends with an $a$. The transition function $p^x$ is

$$p^x(b, a) = p^x(a, b) = 1 - p^x(b, b) = 1 - p^x(a, a) = 1/N.$$

Given a Markov chain with transition function $p^x$ and initial state $a$, the probability that $\mathbf{z}_n = a$, for every $n \leq 2N + 1$, is bounded away from zero. In particular, condition (B1) will not hold. More generally, no homogeneous Markov chain satisfies both (B1) and (B2) in this example.

This example highlights the heart of the problem. The naive solution does satisfy (B1) and (B2) when $p^x$ is sufficiently mixing. However, when it is not, there is no Markov chain that approximates the given sequence in the sense of (B1) and (B2).

The proof of Theorem 1 is organized as follows. First, we provide a general structure result in Section 2.2. When $C$ is a subset of $S$, and $x = (x_1, \ldots, x_N)$ is a sequence of elements in $S$, a $C$-run is a subsequence $(x_{n_1}, x_{n_1+1}, \ldots, x_{n_2})$ such that all its elements are in $C$, while $x_{n_1-1}$ and $x_{n_2+1}$ are not in $C$ (if $n_1 = 1$ or $n_2 = N$, the last condition is vacuous). Our structure result states that given any finite sequence $x$ of elements of $S$, there is a partition of $S$ with the property that for every atom $C$ of the partition and every proper subset $D$ of $C$, the number of $C$-runs is much smaller than the number of $D$-runs. Thus, the sequence moves around inside any atom much more quickly than from one atom to another.

We will use the structure result to argue that the observed transition function $p^x$ associated to $x$, when restricted to any atom of the partition, is mixing. We then construct the simple process $\mathbf{z}$ that approximates $x$. This process will have the following features: (i) it visits every atom $C$ of the partition only *once*, (ii) the duration of the visit to $C$ is $\sum_{s \in C} N_s^x$, the observed number of stages spent in $C$, and (iii) the transition function of $\mathbf{z}$ during the visit to $C$ is $p^x$, properly modified so as to prevent the chain from exiting $C$ too early.

Section 2.3 contains several results on Markov chains. The proof of Theorem 1 is given in Section 2.4.

2.2. *A structure theorem.* We here collect some general notation that is in use throughout the paper. We use the letters $p$ and $q$, with possible sub- or superscripts, to denote transition functions. Probability measures over $S$ are denoted by $\mu$, empirical occupancy measures over $S$ are denoted by $\nu$, while probability measures over $S^{\mathbf{N}}$ are denoted by $\mathbf{P}$. Finally, random variables are often boldfaced, while generic variables are not.

Let a finite sequence $x = (x_0, \ldots, x_N)$ in $S$ be given. For every two subsets $A, B \subseteq S$, we set

$$N_{A,B}^x = \sum_{s \in A, t \in B} N_{s,t}^x \quad \text{and} \quad N_A^x = N_{A,S}^x = \sum_{s \in A} N_s^x.$$

These are the number of transitions from $A$ to $B$ along $x$, and the number of visits to the set $A$ along $x$, respectively. For $C \subseteq S$, we define

$$R_C^x = N_{\overline{C},C}^x + \mathbb{1}_{x_0 \in C}.$$

This is the number of $C$-runs along $x$ [see Feller (1968), II.5]. Plainly, $R_{C \setminus D}^x \leq R_C^x + R_D^x$ for every $D \subset C$, and $|R_C^x - R_{\overline{C}}^x| \leq 1$. Note also that $R_C^x = N_{C,\overline{C}}^x + \mathbb{1}_{x_N \in C}$.

We now state our structure result.

THEOREM 2. *Let $a > 0$ and a finite sequence $x$ of elements of $S$ be given. There is a partition $\mathcal{C}$ of $S$ such that, for every $C \in \mathcal{C}$:*

(P1) $R_C^x \leq (a+1)^{|C|}$.
(P2) *For each $D \subset C$, $R_D^x > aR_C^x$.*

PROOF. Since $R_S^x = 1$, the trivial partition $\mathcal{C} = \{S\}$ satisfies (P1). Among all the partitions that satisfy (P1), let $\mathcal{C}$ be one with maximal number of atoms, and set $k = |\mathcal{C}|$, the number of atoms in $\mathcal{C}$. We will prove that $\mathcal{C}$ satisfies (P2). If it does not, there are $C \in \mathcal{C}$, and a proper subset $D$ of $C$, such that $R_D^x \leq aR_C^x$.

Consider now the partition $\mathcal{C} \setminus \{C\} \cup \{D, C \setminus D\}$ obtained by further partitioning the set $C$ into $D$ and $C \setminus D$. We show that this new partition, with $k+1$ elements, satisfies (P1) as well, contradicting the maximality of $\mathcal{C}$. Indeed, $R_D^x \leq aR_C^x \leq (a+1)^{k+1}$, and $R_{C \setminus D}^x \leq R_C^x + R_D^x \leq R_C^x(a+1) \leq (a+1)^{k+1}$. □

As Theorem 2 has its own merit, we provide two comments concerning the partition $\mathcal{C}$ that satisfies (P1) and (P2).

COMMENT. There need not be a *unique* partition that satisfies both (P1) and (P2). Indeed, let $S = \{0, 1\}$ and $x = (0, 1, 0, 1, \ldots, 0, 1)$ (a sequence of length $N+1$), and let $a > 0$ be such that $a < \frac{N+1}{2} \leq (a+1)^2$. Since $R_{\{0\}}^x = R_{\{1\}}^x = \frac{N+1}{2}$, the two partitions of $S$ satisfy (P1) and (P2).

COMMENT. For $a \geq 2$, the partition that is defined in the proof of Theorem 2 is unique. To verify this claim, it is enough to check that, given two partitions $\mathcal{C}$ and $\mathcal{D}$ that satisfy (P2), the following holds: for every $C \in \mathcal{C}$ and $D \in \mathcal{D}$, if the intersection $C \cap D$ is nonempty, then it is equal to either $C$ or $D$.

Assume to the contrary that $P = C \cap D$ is a proper subset of both $C$ and $D$. Then one has $N_{P,\overline{D}}^x + N_{P,D \setminus P}^x = N_{P,\overline{C}}^x + N_{P,C \setminus P}^x = R_P^x - \mathbb{1}_{x_N \in P}$, $N_{P,C \setminus P}^x + N_{P,D \setminus P}^x \leq R_P^x - \mathbb{1}_{x_N \in P}$, $N_{P,\overline{C}}^x \leq N_{C,\overline{C}}^x = R_C^x - \mathbb{1}_{x_N \in C}$, and $N_{P,\overline{D}}^x \leq N_{D,\overline{D}}^x = R_D^x - \mathbb{1}_{x_N \in D}$. It follows that

$$R_P^x - \mathbb{1}_{x_N \in P} \geq N_{P,C \setminus P}^x + N_{P,D \setminus P}^x$$
$$= 2R_P^x - 2 \times \mathbb{1}_{x_N \in P} - N_{P,\overline{C}}^x - N_{P,\overline{D}}^x$$
$$\geq 2R_P^x - R_C^x - R_D^x.$$

In particular, by (P2),

$$R_C^x + R_D^x - \mathbb{1}_{x_N \in P} \geq R_P^x > a \times \max\{R_C^x, R_D^x\},$$

a contradiction when $a \geq 2$.

2.3. *On Markov chains.* We here collect a few useful results about Markov chains. First, we provide a result on the speed of convergence of an irreducible Markov chain to its invariant measure. Next, we make a few observations on the expected exit time from sub-domains of $S$.

Throughout the present section, a transition function $p$ over $S$ is given. For $s \in S$, we denote by $\mathbf{P}_{s,p}$ the law of a homogeneous Markov chain $\mathbf{z}$ with transition function $p$ and initial state $s$, and by $\mathbf{E}_{s,p}$ the expectation w.r.t. $\mathbf{P}_{s,p}$. For $\mu \in \mathcal{P}(S)$, $\mathbf{E}_{\mu,p} = \sum_{s \in S} \mu(s) \mathbf{E}_{s,p}$ is the expectation operator when the initial state is chosen according to $\mu$.

The *hitting time* of a set $C \subseteq S$ is $T_C = \min\{n \geq 0 : \mathbf{z}_n \in C\}$ (with $\min \varnothing = +\infty$). For $t \in S$, we abbreviate $T_{\{t\}}$ to $T_t$ and we denote by $T_t^+ = \min\{n \geq 1 : \mathbf{z}_n = t\}$ the *first return time* to $t$.

2.3.1. *Convergence to the invariant measure.*

DEFINITION 2. Let $\gamma > 0$ be given. The transition function $p$ is $\gamma$-*mixing* if $\mathbf{E}_{s,p}[T_t^+] \leq \gamma$, for every $s, t \in S$.

Plainly, a $\gamma$-mixing transition function is irreducible. The next theorem bounds the speed of convergence of the empirical occupation measure to the invariant measure for $\gamma$-mixing homogeneous Markov chains. In this statement $\bar{\nu}_m^{\mathbf{z}}$ is the occupancy measure in stages 1 through $m$: $\bar{\nu}_m^{\mathbf{z}}(s) = \frac{1}{m}|\{1 \leq n \leq m : \mathbf{z}_n = s\}|$.

THEOREM 3. *Assume that the transition function $p$ is $\gamma$-mixing and let $\mu$ be its invariant measure. Let $n \in \mathbf{N}$ and $\varepsilon \in (0, 1/2)$ be such that $\varepsilon n > 4\gamma$. Then, for every $s, t \in S$,*

$$\mathbf{P}_{t,p}\left(\left|\frac{\bar{\nu}_n^{\mathbf{z}}(s)}{\mu(s)} - 1\right| \geq \varepsilon\right) < \frac{17\gamma}{n\varepsilon^2}. \tag{2}$$

REMARK 2. Inspection of the proof shows that inequality (2) holds more generally for each state $s \in S$ such that $\max_{t \in S} \mathbf{E}_{t,p}[T_s^+] \leq \gamma$.

REMARK 3. Since $|\bar{\nu}_n^{\mathbf{z}}(s) - \nu_n^{\mathbf{z}}(s)| \leq \frac{1}{n}$, one has, under the assumptions of Theorem 3,

$$\mathbf{P}_{t,p}\left(|\nu_n^{\mathbf{z}}(s) - \mu(s)| \geq \varepsilon\mu(s) + \frac{1}{n}\right) < \frac{17\gamma}{n\varepsilon^2}. \tag{3}$$

REMARK 4. It is likely that the bound in Theorem 3 can be substantially improved, possibly to an exponential bound. Recently, Glynn and Ormoneit (2002) provided a generalization of Hoeffding's inequality to uniformly ergodic chains. However, their ergodicity assumption (A1) is stronger than our mixing assumption, hence our result does not follow from their statement.

PROOF OF THEOREM 3. The proof relies on the following two identities:

$$(4) \quad \mathbf{E}_{s,p}[T_s^+] = \frac{1}{\mu(s)} \quad \text{and} \quad \mu(s)\text{Var}_{s,p}(T_s^+) = 2\mathbf{E}_{\mu,p}[T_s] + 1 - \frac{1}{\mu(s)}$$

[see Aldous and Fill (2002), Chapter 2, identity (22) for the second one]. Since $p$ is $\gamma$-mixing, $1/\mu(s) = \mathbf{E}_{s,p}[T_s^+] \leq \gamma < \varepsilon n/4$, and $\mathbf{E}_{\mu,p}[T_s] \leq \mathbf{E}_{\mu,p}[T_s^+] \leq \gamma$. Since $1 - 1/\mu(s) \leq 0$, we also have $\mu(s)\text{Var}_{s,p}(T_s^+) \leq 2\gamma$.

For notational clarity, we set $n_\varepsilon = \lceil n\mu(s)(1-\varepsilon) \rceil$ and $n^\varepsilon = \lfloor n\mu(s)(1+\varepsilon) \rfloor$. Note that $n_\varepsilon \leq n^\varepsilon$. Moreover, $n_\varepsilon + n^\varepsilon - 1 \leq 2n\mu(s)$, so that $n_\varepsilon + n^\varepsilon \leq 3n\mu(s)$.

On $\{\bar{\nu}_n^\mathbf{z}(s) \leq \mu(s)(1-\varepsilon)\}$ one has $V_{s,n_\varepsilon} \geq n$, whereas on $\{\bar{\nu}_n^\mathbf{z}(s) \geq \mu(s)(1+\varepsilon)\}$ one has $V_{s,n^\varepsilon} \leq n$. Therefore, the event $\{|\frac{\bar{\nu}_n^\mathbf{z}(s)}{\mu(s)} - 1| \geq \varepsilon\}$ is included in the union of the two events $\{V_{s,n_\varepsilon} \geq n\}$ and $\{V_{s,n^\varepsilon} \leq n\}$, so that

$$(5) \quad \mathbf{P}_{t,p}\left(\left|\frac{\bar{\nu}_n^\mathbf{z}(s)}{\mu(s)} - 1\right| \geq \varepsilon\right) \leq \mathbf{P}_{t,p}(V_{s,n_\varepsilon} \geq n) + \mathbf{P}_{t,p}(V_{s,n^\varepsilon} \leq n).$$

We will prove the result by providing an upper bound on the probability that $V_{s,n_\varepsilon} \geq n$ and on the probability that $V_{s,n^\varepsilon} \leq n$.

Since $\varepsilon < 1/2$ and since $\varepsilon n \mu(s) > 4$, straightforward manipulations show that

$$(6) \quad \min\left\{n - \frac{n_\varepsilon}{\mu(s)}, \frac{n^\varepsilon}{\mu(s)} - n\right\} \geq \frac{3}{4}n\varepsilon \quad \text{and}$$

$$(7) \quad \min\left\{n(1-\varepsilon^2) - \frac{n_\varepsilon - 1}{\mu(s)}, \frac{n^\varepsilon - 1}{\mu(s)} + 1 - n\right\} \geq \frac{1}{2}n\varepsilon.$$

For $s \in S$ and $k \in \mathbf{N}$, let $V_{s,k}$ denote the time of the $k$th return to $s$ (with $V_{s,0} = 0$), and let $T_{s,k}^+ = V_{s,k} - V_{s,k-1}$ denote the length of the $k$th visit to $S \setminus \{s\}$.

We distinguish the two cases $s = t$ and $s \neq t$.

CASE 1. $s = t$.

In this case the random variables $T_{s,k}^+$ are i.i.d. and share the law of $T_s^+$. Since $\mathbf{E}_{s,p}[T_s^+] = 1/\mu(s)$, one has, by Chebyshev's inequality and by (6),

$$\mathbf{P}_{s,p}(V_{s,n_\varepsilon} \geq n) = \mathbf{P}_{s,p}\left(V_{s,n_\varepsilon} - \frac{n_\varepsilon}{\mu(s)} \geq n - \frac{n_\varepsilon}{\mu(s)}\right) \leq \frac{n_\varepsilon \text{Var}_{s,p}(T_s^+)}{(3/4n\varepsilon)^2}$$

and

$$\mathbf{P}_{s,p}(V_{s,n^\varepsilon} \leq n) \leq \mathbf{P}_{s,p}\left(\frac{n^\varepsilon}{\mu(s)} - V_{s,n^\varepsilon} \geq \frac{n^\varepsilon}{\mu(s)} - n\right) \leq \frac{n^\varepsilon \text{Var}_{s,p}(T_s^+)}{(3/4n\varepsilon)^2}.$$

Hence, by (5), (4) and since $p$ is $\gamma$-mixing, one obtains

$$(8) \quad \mathbf{P}_{s,p}\left(\left|\frac{\bar{\nu}_n^\mathbf{z}(s)}{\mu(s)} - 1\right| \geq \varepsilon\right) \leq \frac{(n_\varepsilon + n^\varepsilon)\text{Var}_{s,p}(T_s^+)}{(3/4n\varepsilon)^2} \leq \frac{16 \times 3 \times 2\gamma}{9n\varepsilon^2} < \frac{17\gamma}{n\varepsilon^2}.$$

CASE 2. $s \neq t$.

By Markov's inequality and since $p$ is $\gamma$-mixing,

$$\mathbf{P}_{t,p}(T_s^+ \geq \varepsilon^2 n) \leq \frac{\gamma}{n\varepsilon^2}. \tag{9}$$

By repeating the steps of Case 1 using (7), one has

$$\begin{aligned}
\mathbf{P}_{t,p}(T_{s,1}^+ &\leq \varepsilon^2 n, V_{s,n_\varepsilon} \geq n) \\
&\leq \mathbf{P}_{s,p}(T_{s,2}^+ + \cdots + T_{s,n_\varepsilon}^+ \geq n(1-\varepsilon^2)) \\
&\leq \frac{(n_\varepsilon - 1)\mathrm{Var}_{s,p}(T_s^+)}{(n(1-\varepsilon^2) - (n_\varepsilon - 1)/\mu(s))^2} \leq \frac{(n_\varepsilon - 1)\mathrm{Var}_{s,p}(T_s^+)}{(1/2n\varepsilon)^2},
\end{aligned} \tag{10}$$

while

$$\begin{aligned}
\mathbf{P}_{t,p}(V_{s,n^\varepsilon} \leq n) &\leq \mathbf{P}_{s,p}(T_{s,2}^+ + \cdots + T_{s,n^\varepsilon}^+ \leq n-1) \\
&\leq \frac{(n^\varepsilon - 1)\mathrm{Var}_{s,p}(T_s^+)}{(1/2n\varepsilon)^2}.
\end{aligned} \tag{11}$$

By summing (9)–(11), one obtains

$$\mathbf{P}_{t,p}\left(\left|\frac{\bar\nu_n^{\mathbf{z}}(s)}{\mu(s)} - 1\right| \geq \varepsilon\right) \leq \frac{4\mathrm{Var}_{s,p}(T_s^+)(n_\varepsilon + n^\varepsilon - 2) + \gamma n}{n^2 \varepsilon^2}.$$

Therefore,

$$\mathbf{P}_{t,p}\left(\left|\frac{\bar\nu_n^{\mathbf{z}}(s)}{\mu(s)} - 1\right| \geq \varepsilon\right) \leq \frac{4 \times 2 \times 2\gamma + \gamma}{n\varepsilon^2} = \frac{17\gamma}{n\varepsilon^2},$$

as desired. □

2.3.2. *Expected exit times.* We here analyze the exit time from a given sub-domain. Our estimates use two new mixing measures for irreducible transition functions.

Throughout this section we assume that the transition function $p$ is irreducible with invariant measure $\mu$. We use repeatedly the inequality

$$\mathbf{E}_{s,p}[T_{\overline{L}}] \leq \mathbf{E}_{s,p}[T_{\overline{L \cup \{t\}}}] + \mathbf{E}_{t,p}[T_{\overline{L}}], \tag{12}$$

which holds for every $L \subset S$ and every $s, t \in L$.

DEFINITION 3. For $C \subseteq S$, we define

$$\lambda_p(C) = \max_{s \in C} \mathbf{E}_{s,p}[T_{\overline{C}}] \quad \text{and} \quad \rho_p(C) = \max_{D \subset C} \min_{s \in D} \mathbf{E}_{s,p}[T_{\overline{D}}].$$

Observe that one always has $\lambda_p(C) \geq \rho_p(C)$. $\lambda_p(C)$ bounds the time it takes to leave $C$. $\rho_p(C)$ may be interpreted as a measure of how fast a Markov chain with transition function $p$ visits each and every state of $C$. The following lemma adds substance to these interpretations. We shall use it to derive further estimates of exit times.

LEMMA 1. *For every $C \subseteq S$:*

(i) $\mathbf{E}_{s,p}[T_{\overline{D}}] \leq |D|\rho_p(C)$, *for every $D \subset C$ and every $s \in D$.*
(ii) $\mathbf{E}_{s,p}[T_{\overline{C}}] \geq \lambda_p(C) - (|C|-1)\rho_p(C)$, *for every $s \in C$.*

PROOF. We prove the first statement by induction over $|D|$. Plainly, the inequality holds for singletons. Assume it holds for every subset with $k$ elements. Let $D \subset C$ be any subset with $k+1$ elements, and let $s \in D$. Choose $t \in D$ such that $\mathbf{E}_{t,p}[T_{\overline{D}}] \leq \rho_p(C)$. By (12) and the induction hypothesis, applied to $D \setminus \{t\}$, we have

$$\mathbf{E}_{s,p}[T_{\overline{D}}] \leq \mathbf{E}_{s,p}[T_{\overline{D} \cup \{t\}}] + \mathbf{E}_{t,p}[T_{\overline{D}}] \leq (|D|-1)\rho_p(C) + \rho_p(C).$$

We now prove the second statement. Let $s \in C$ be given, and let $t \in C$ be such that $\lambda_p(C) = \mathbf{E}_{t,p}[T_{\overline{C}}]$. If $t = s$, (ii) trivially holds. Otherwise, by (12) and (i),

$$\mathbf{E}_{s,p}[T_{\overline{C}}] \geq \mathbf{E}_{t,p}[T_{\overline{C}}] - \mathbf{E}_{t,p}[T_{\overline{C} \cup \{s\}}] \geq \lambda_p(C) - (|C|-1)\rho_p(C),$$

as desired. □

The following lemma bounds the probability that the process leaves a set $C$ before it visits some given state $t \in C$.

LEMMA 2. *For every $C \subset S$ and every $s, t \in C$, one has*

$$(13) \qquad \mathbf{P}_{s,p}(T_{\overline{C}} < T_t) \leq 2|C|\frac{\rho_p(C)}{\lambda_p(C) - (|C|-1)\rho_p(C)}.$$

PROOF. If $s = t$, the left-hand side in (13) vanishes, so that the lemma trivially holds. Hence we assume from now on that $s \neq t$, so that $|C| \geq 2$, and, therefore, $\rho_p(C) \geq 1$.

We modify the state space $S$ and the transition function $p$ by collapsing $\overline{C}$ to a single state, still denoted $\overline{C}$, which leads to $t$ in one step. Since this change does not affect the probability that $T_{\overline{C}} < T_t$, we still denote the modified transition function by $p$. This amounts to assuming that $p(\overline{C}, t) = 1$, hence $\mathbf{E}_{\overline{C},p}[T_t] = 1$. By Aldous and Fill [(2002), Chapter 2, Corollary 10],

$$(14) \qquad \mathbf{P}_{s,p}(T_{\overline{C}} < T_t) = \frac{\mathbf{E}_{s,p}[T_t] + \mathbf{E}_{t,p}[T_{\overline{C}}] - \mathbf{E}_{s,p}[T_{\overline{C}}]}{\mathbf{E}_{\overline{C},p}[T_t] + \mathbf{E}_{t,p}[T_{\overline{C}}]}.$$

Since $\mathbf{E}_{\overline{C},p}[T_t] = 1$, one has, by (12), $\mathbf{E}_{s,p}[T_t] \leq \mathbf{E}_{s,p}[T_{\overline{C}\cup\{t\}}] + 1$. Equation (12) also implies that $\mathbf{E}_{t,p}[T_{\overline{C}}] - \mathbf{E}_{s,p}[T_{\overline{C}}] \leq \mathbf{E}_{t,p}[T_{\overline{C}\cup\{s\}}]$. Therefore, by Lemma 1(i), the numerator in (14) is at most

$$1 + \mathbf{E}_{t,p}[T_{\overline{C}\cup\{s\}}] + \mathbf{E}_{s,p}[T_{\overline{C}\cup\{t\}}] \leq 2(|C|-1)\rho_p(C) + 1 < 2|C|\rho_p(C).$$

On the other hand, the denominator is equal to $1 + \mathbf{E}_{t,p}[T_{\overline{C}}]$, hence by Lemma 1(ii) is at least $\lambda_p(C) - (|C|-1)\rho_p(C)$. $\square$

For $C \subset S$, we denote by $p_C$ the transition function $p$ watched on $C$ [see Aldous and Fill (2002), Chapter 2, Section 7.1]. Formally,

$$(15) \quad p_C(s,t) = p(s,t) + \sum_{u \notin C} p(s,u) \mathbf{P}_{u,p}(T_C = T_t) \qquad \text{for every } s,t \in C.$$

Since $p$ is irreducible, the transition function $p_C$ is irreducible, and its invariant measure coincides with the invariant measure $\mu$ of $p$, conditioned on $C$: $\mu(s|C) = \mu(s)/\mu(C)$, for every $s \in C$ [see Aldous and Fill (2002)].

The next lemma bounds the time it takes the process to reach a given state $t \in C \subseteq S$, when watched on $C$. Thus, it bounds the expected number of stages the Markov chain with transition function $p$ spends in $C$ until it reaches $t$ for the first time.

LEMMA 3. *For $s,t \in C$, one has $\mathbf{E}_{s,p_C}[T_t] \leq \frac{(|C|-1)\rho_p(C)}{1-\max_{u\in C}\mathbf{P}_{u,p}(T_{\overline{C}}<T_t)}$.*

PROOF. If $s = t$ the lemma trivially holds, as in this case $\mathbf{E}_{s,p_C}[T_t] = 0$.

Assume then that $s \neq t$, so, in particular, $|C| \geq 2$. Let $t \in S$ be given. For convenience set $\alpha = \max_{s \in C} \mathbf{E}_{s,p_C}[T_t]$, and let $s' \in S$ achieve the maximum. Since $|C| \geq 2$, $s' \neq t$. Therefore, by Lemma 1(i),

$$\alpha = \mathbf{E}_{s',p_C}[T_t] \leq \mathbf{E}_{s',p}[T_{\overline{C}\cup\{t\}}] + \mathbf{P}_{s',p}(T_{\overline{C}} < T_t)\alpha$$
$$\leq (|C|-1)\rho_p(C) + \alpha \mathbf{P}_{s',p}(T_{\overline{C}} < T_t).$$

Thus, for every $s \in C$,

$$\mathbf{E}_{s,p_C}[T_t] \leq \alpha \leq \frac{(|C|-1)\rho_p(C)}{1 - \mathbf{P}_{s',p}(T_{\overline{C}} < T_t)} \leq \frac{(|C|-1)\rho_p(C)}{1 - \max_{u \in C} \mathbf{P}_{u,p}(T_{\overline{C}} < T_t)},$$

as desired. $\square$

We conclude with two results, stated without proof. First, given $C \subset S$, define

$$(16) \quad K_p(C) = \frac{\sum_{s \in C} \mu(s)}{\sum_{s \in C} \mu(s) p(s, \overline{C})}.$$

The numerator in (16) is the frequency of stages spent in $C$, while the denominator is the frequency of exits from $C$. Therefore, $K_p(C)$ is the average length of a visit to $C$. In particular, the following holds:

$$(17) \qquad \min_{s \in C} \mathbf{E}_{s,p}[T_{\overline{C}}] \leq K_p(C) \leq \max_{s \in C} \mathbf{E}_{s,p}[T_{\overline{C}}] = \lambda_p(C).$$

Second, straightforward computations show that for every $C \subset S$,

$$(18) \qquad \frac{\mu(C)}{\mu(\overline{C})} \geq \frac{\min_{s \in C} \mathbf{E}_{s,p}[T_{\overline{C}}]}{\max_{s \in \overline{C}} \mathbf{E}_{s,p}[T_C]}.$$

2.4. *Proof of Theorem* 1. This section is devoted to the proof of Theorem 1. It is convenient to deal with sequences $x$ that are *exhaustive* ($N_s^x > 0$ for each $s \in S$) and *periodic* ($x_N = x_0$). The general result will follow since an arbitrary sequence $x$ can be extended into an exhaustive and periodic one, by appending at most $|S|$ elements to $x$ (see details in Section 2.4.3). Observe that for the purpose of Theorem 1, all sequences can be assumed to be exhaustive, since states that are not visited along $x$ can simply be dropped. Since this assumption cannot be made to prove the more general theorem of this paper, we prefer not to make it here as well.

The assumption that the sequence is exhaustive and periodic allows us to make use of the following lemma whose proof is omitted.

LEMMA 4. *Let $x = (x_0, \ldots, x_N)$ be exhaustive and periodic. The observed transition function $p^x$ is irreducible, and its invariant measure coincides with the observed occupancy measure $\nu^x(s) = \frac{N_s^x}{N}$.*

Let $\varepsilon \in (0, \frac{1}{2})$, $\delta \in (0, \frac{1}{4|S|+1})$ and $\zeta \in (0, 2\delta)$ be given. We choose $N_0 \in \mathbf{N}$ large enough so that (N1) $N_0^{1-(4|S|+1)\delta} > 2^{|S|}/\varepsilon$, and (N2) $N_0^{2\delta-\zeta} > 4 \times 17|S|^4/\varepsilon^2$. Therefore, we have, in particular, (N3) $(N_0^{4\delta}+1)^{|S|}N_0^{\delta-1} \leq \varepsilon/(2^{|S|}+1)$, (N4) $N_0^{\delta} > \max\{|S|^2/\varepsilon, 2|S|+2\}$, (N5) $N_0^{2\delta} \geq 8/\varepsilon$ and (N6) $N_0^{4\delta} \geq \max\{1/\varepsilon, 10|S|\}$.

We will prove that the conclusion of Theorem 1 holds for every $N \geq N_0$ and every exhaustive, periodic sequence $x$. We first apply Theorem 2 to the sequence $x$, with $a = N^{4\delta}$, to obtain a partition $\mathcal{C} = (S_1, \ldots, S_K)$ of $S$ that satisfies the conclusions of that theorem. Observe that $a$ depends on the length of the sequence. We now proceed as follows. In Section 2.4.1 we argue that the transition function $p^x$ is mixing, when watched on any atom $S_k$ of $\mathcal{C}$. In Section 2.4.2 we define the approximating process, and we check that assertions (B1) and (B2) hold.

2.4.1. *Properties of $p_{S_k}^x$.* Following the notation in use in Section 2.3.2, we denote by $p_{S_k}^x$ the transition function $p^x$, when watched on $S_k$. Since $p^x$ is irreducible, so is $p_{S_k}^x$. The goal of this section is to prove that $p_{S_k}^x$ is $N^{1-3\delta}$-mixing (see Proposition 1 below). To this end, we first relate the mixing constants $\lambda_{p^x}(S_k)$ and $\rho_{p^x}(S_k)$ to the features of $x$.

LEMMA 5. *Let $k$ be such that $|S_k| > 1$. One has*

$$\rho_{p^x}(S_k) \leq \max_{D \subset S_k} \frac{N_D^x}{R_D^x - 1} \leq \frac{2}{a} \lambda_{p^x}(S_k). \tag{19}$$

Note that the quantity $\frac{N_D^x}{R_D^x - 1}$ is approximately the average length of a visit to $D$ along $x$. Thus, the expected exit time from $D \subset S_k$ is much smaller than the expected exit time from $S_k$.

PROOF OF LEMMA 5. Let $k$ be such that $|S_k| > 1$, and let $D \subseteq S_k$ (here $D$ may be equal to $S_k$). By Lemma 4 one has $\nu^x(D) = N_D^x/N$ and $\sum_{s \in D} \nu^x(s) p^x(s, \overline{D}) = N_{D,\overline{D}}^x/N$. By (16), as long as $D \subset S$, $K_{p^x}(D) = N_D^x/N_{D,\overline{D}}^x$, so that

$$\frac{N_D^x}{R_D^x} \leq K_{p^x}(D) \leq \frac{N_D^x}{R_D^x - 1}. \tag{20}$$

If $D$ is a strict subset of $S_k$, (17), the second inequality in (20), (P2) and (N4) yield

$$\min_{s \in D} \mathbf{E}_{s,p^x}[T_{\overline{D}}] \leq K_{p^x}(D) \leq \frac{N_D^x}{R_D^x - 1} \leq \frac{N_{S_k}^x}{aR_{S_k}^x - 1} \leq \frac{2}{a} \times \frac{N_{S_k}^x}{R_{S_k}^x}. \tag{21}$$

The left-hand side inequality in (19) follows by taking the maximum over $D \subset S_k$.

If $S_k = S$, $\lambda_{p^x}(S_k) = +\infty$, and the right-hand side inequality in (19) trivially holds. Otherwise, when applied to $S_k$, the first inequality in (20) and (17) yield

$$\frac{N_{S_k}^x}{R_{S_k}^x} \leq K_{p^x}(S_k) \leq \lambda_{p^x}(S_k), \tag{22}$$

and the right-hand side inequality in (19) follows from (21) and (22). □

PROPOSITION 1. *If $|S_k| \geq 2$, the transition function $p_{S_k}^x$ is $N^{1-3\delta}$-mixing.*

PROOF. We will prove that $\mathbf{E}_{s,p_{S_k}^x}[T_t] \leq N^{1-3\delta} - 1$, for each $s, t \in S_k$. Let $s, t \in S_k$ be given. By Lemma 3,

$$\mathbf{E}_{s,p_{S_k}^x}[T_t] \leq \frac{(|S_k| - 1)\rho_{p^x}(S_k)}{1 - \max_{u \in S_k} \mathbf{P}_{u,p^x}(T_{\overline{S_k}} < T_t)}.$$

By Lemma 2, the denominator is at least $1 - 2|S_k|\frac{\rho_{p^x}(S_k)}{\lambda_{p^x}(S_k) - (|S_k| - 1)\rho_{p^x}(S_k)} \geq \frac{1}{2}$, where the inequality holds by Lemma 5 and (N6). Therefore,

$$\mathbf{E}_{s,p^x_{S_k}}[T_t] \leq 2|S_k|\rho_{p^x}(S_k). \tag{23}$$

By Lemma 5, (P2) and (N4),

$$\rho_{p^x}(S_k) \leq \max_{D \subset S_k} \frac{N^x_D}{R^x_D - 1} < \frac{N}{N^{4\delta} - 1} \leq \frac{N^{1-3\delta} - 1}{2|S_k|}. \tag{24}$$

The result follows by combining (23) and (24). $\square$

2.4.2. *The approximating process.* We now construct a Markov chain $\mathbf{z}$ that approximates the sequence $x$. Ideally the chain is composed of $|\mathcal{C}|$ pieces, with the length of piece $k$ being $N^x_{S_k}$. However, to avoid degenerate cases, we take into account only the atoms $S_k$ that are frequently visited, that is, those with $N^x_{S_k} \geq N^{1-\delta}$.

Set $K_0 = \{k : N^x_{S_k} \geq N^{1-\delta}\}$, and assume for convenience that $K_0$ contains the first $|K_0|$ atoms in $\mathcal{C}$, so that $K_0 = \{1, \ldots, |K_0|\}$. Assume, moreover, that the set $S_{|K_0|}$ is the most frequently visited set, so that, in particular, $N^x_{S_{|K_0|}} \geq N/|S|$.

The chain $\mathbf{z}$ is a piecewise homogeneous Markov chain with $|K_0|$ pieces. The "extra" stages that are created by the removal of rarely visited atoms in $\mathcal{C}$ are added to piece $|K_0|$. Since piece $|K_0|$ is the most frequently visited piece, this will hardly affect the estimates for that piece.

Formally, we denote by $m_k$ the length of piece $k$. Thus, $m_k = N^x_{S_k}$ if $k < |K_0|$, and $m_k = N - \sum_{j < |K_0|} m_j$ if $k = |K_0|$. In particular, $1 \leq \frac{m_k}{N^x_{S_k}} \leq 1 + \frac{|S|^2}{N^\delta}$. For $k = 1, \ldots, K_0$, we let $p_k$ be a transition function such that

$$\begin{aligned} p_k(s,t) &= p^x_{S_k}(s,t), & s \in S_k, t \in S_k, \\ p_k(s, S_k) &= 1, & s \notin S_k. \end{aligned}$$

The exact definition of $p_k(s, \cdot)$ for $s \notin S_k$ is irrelevant. Thus, $p_k$ moves in one step to $S_k$ and coincides with $p^x_{S_k}$ there.

We let $\mathbf{z}$ be a Markov chain with initial state in $S_1$, and transitions $p^{\mathbf{z}}_n(s,t) = p_k(s,t)$ if $-1 + \sum_{j<k} m_j \leq n < -1 + \sum_{j \leq k} m_j$. (This does not define $p^{\mathbf{z}}_{N-1}$. The choice of $p^{\mathbf{z}}_{N-1}$ is irrelevant. On the other hand, it defines $p^{\mathbf{z}}_{-1}$, which is never used.) Note that $\mathbf{z}$ is a piecewise homogeneous Markov chain with $|K_0| \leq |S|$ pieces.

Plainly, the chain $\mathbf{z}$ visits each set $S_k$, $k \leq K_0$, only once, for exactly $m_k$ stages: from stage $\sum_{j<k} m_j$ to stage $\sum_{j \leq k} m_j - 1$ (inclusive).

We now prove that both assertions in Theorem 1 hold.

We start with assertion (B1). Let $k \in K$ and $s \in S_k$ be given. We discuss three cases.

If $N^x_{S_k} < N^{1-\delta}$, then $\nu^x(s) = N^x_s/N < 1/N^\delta$, and (B1) trivially holds.

Assume now that $N^x_{S_k} \geq N^{1-\delta}$ and that $S_k = \{s\}$ is a singleton. By construction, the chain $\mathbf{z}$ will be in state $s$ all through piece $k$, and only in those stages. In particular, by (N4),

$$0 \leq \frac{\nu^{\mathbf{z}}_N(s)}{\nu^x(s)} - 1 \leq \frac{m_k - N^x_s}{N^x_s} \leq \frac{|S|^2}{N^\delta} < \varepsilon,$$

and (B1) holds.

Assume finally that $N^x_{S_k} \geq N^{1-\delta}$ and $|S_k| \geq 2$. By (N5) and Proposition 1, the assumptions of Theorem 3 hold w.r.t. $p = p^x_{S_k}$, $\gamma = N^{1-3\delta}$, $n = m_k \geq N^{1-\delta}$ and $\varepsilon/2$. For every $t \in S_k$, one has, by (N5), Remark 3 and (N2),

$$\begin{aligned}
\mathbf{P}_{t,p_k}\left(\left|\frac{\nu^{\mathbf{z}}_{m_k}(s)}{\nu^x(s|S_k)} - 1\right| \geq \varepsilon\right) \\
(25) \quad \leq \mathbf{P}_{t,p_k}\left(\left|\frac{\nu^{\mathbf{z}}_{m_k}(s)}{\nu^x(s|S_k)} - 1\right| \geq \frac{\varepsilon}{2} + \frac{1}{m_k\nu^x(s|S_k)}\right) \\
\leq 4 \times \frac{17N^{1-3\delta}}{m_k\varepsilon^2} \leq 4 \times \frac{17N^{1-3\delta}}{\varepsilon^2 N^{1-\delta}} \leq \frac{1}{|S|} \times \frac{1}{N^\zeta}.
\end{aligned}$$

Since the chain $\mathbf{z}$ does not visit $s \in S_k$, except in piece $k$, (B1) follows from (25).

We now turn to assertion (B2). Since $\mathbf{z}$ never visits states in $\bigcup_{k \notin K_0} S_k$, we need to verify that (B2) holds for states $s \in \bigcup_{k \in K_0} S_k$. Observe that for every $k \in K_0$ and every $s \in S_k$, one has

$$\|p_k(s,\cdot) - p^x(s,\cdot)\| \leq \sum_{u \notin S_k} p^x(s,u) = \frac{N^x_{s,\overline{S_k}}}{N^x_s} \leq \frac{R^x_{S_k}}{N^x_s}.$$

If $S_k = \{s\}$ is a singleton, then by Theorem 2(P1) and (N3), the right-hand side is bounded by $\frac{(N^{4\delta}+1)^{|S|}}{N^{1-\delta}} < \varepsilon$. If $|S_k| \geq 2$, then since $N^x_s \geq R^x_s$, by Theorem 2(P2) and (N6), the right-hand side is bounded by $\frac{R^x_{S_k}}{R^x_s} \leq \frac{R^x_{S_k}}{N^{4\delta} R^x_{S_k}} < \varepsilon$.

Since for every stage $n$ in piece $k$, except the last one, $p^{\mathbf{z}}_n = p_k$ and $\mathbf{z}_n \in S_k$, $\mathbf{P}_{\mathbf{z}}$-a.s., one has for such $n$'s,

$$\|p^{\mathbf{z}}_n(\mathbf{z}_n,\cdot) - p^x(\mathbf{z}_n,\cdot)\| = \|p_k(\mathbf{z}_n,\cdot) - p^x(\mathbf{z}_n,\cdot)\| \leq \varepsilon \quad \text{almost surely.}$$

2.4.3. *The case of arbitrary sequences.* Our goal now is to prove Theorem 1 for any sequence of observations. We add a few stages to the sequence $x$—fictitious observations — in order to obtain a periodic and exhaustive sequence $x^*$. We then apply the above analysis to the augmented sequence $x^*$. Since only few fictitious observations are needed to change $x$ into a periodic and exhaustive sequence, the desired result will follow.

We let $\varepsilon \in (0, 1/2)$, $\delta \in (0, \frac{1}{2(4|S|+1)})$ and $\zeta \in (0, 2\delta)$ be given. Set $\delta' = 2\delta \in (0, \frac{1}{4|S|+1})$ and $\varepsilon' = \varepsilon - 2\varepsilon^2$.

Choose $N_0 \in \mathbf{N}$ such that (N1) and (N2) hold for $N_0$ w.r.t. $\delta'$, $\varepsilon'$ and $\zeta$ (rather than w.r.t. $\delta$, $\varepsilon$ and $\zeta$). We will argue that the conclusion of Theorem 1 holds for every $N \geq N_0$.

Let $x = (x_0, \ldots, x_N)$ be an arbitrary sequence in $S$, and let $S^* = \bigcup_{n=0}^{N} \{x_n\} \subseteq S$ be the set of states visited by $x$. Consider the sequence $x^* = (x_0, x_1, \ldots, x_N, x_1^*, \ldots, x_r^*, x_0)$, where $r = |S| - |S^*|$ is the number of states not visited by $x$, and $S \setminus S^* = \{x_1^*, \ldots, x_r^*\}$. By construction, this new sequence is periodic and exhaustive. The length $N_* + 1$ of this sequence is $N + r + 2 < N + |S| + 2$.

By construction $N_*$ satisfies (N1) and (N2) with $\delta'$ and $\varepsilon'$. Therefore, there is a piecewise Markov chain $(\mathbf{z}_n)_{n \leq N_*}$ such that (B1) and (B2) hold w.r.t. $N_*$, $\nu^{x^*}$, $\varepsilon'$, $\delta'$ and $\zeta$. Observe that each state $x_j^* \in S \setminus S^*$ constitutes a singleton in the partition $\mathcal{C}$ associated with $x^*$, and that $N_{x_j^*}^{x^*} = 1$, so that it is never visited by $\mathbf{z}$.

One can verify that the restriction of $\mathbf{z}$ to the first $N$ stages satisfies (B1) and (B2) w.r.t. $N$, $\nu^x$, $\varepsilon$, $\delta$ and $\zeta$. The computations are tedious and of no specific interest, and are therefore omitted.

## 3. The general problem.

3.1. *Presentation and discussion.* We here address the more general problem of devising an approximating simple process, given structural constraints on the process. In other words, we wish to construct a simple process within a given class of processes. The kind of structural constraints we allow for is described as follows. For each $s \in S$, we let a nonempty polyhedron $V(s) \subseteq \mathcal{P}(S)$ be given. Recall that a *polyhedron* is the convex hull of finitely many points. Let $V = (V(s))_{s \in S}$ denote the product polyhedron, and for every $s \in S$ denote by $V^*(s)$ the set of extreme points of $V(s)$.

DEFINITION 4. A *V-process* is an $S$-valued process $\mathbf{z} = (\mathbf{z}_n)_n$ such that for every $n \geq 0$ the conditional distribution of $\mathbf{z}_n$, given $\mathbf{z}_0, \ldots, \mathbf{z}_{n-1}$, is in $V(\mathbf{z}_{n-1})$, $\mathbf{P}_\mathbf{z}$-a.s.

In a sense, one-step transitions are required to satisfy exogeneously given constraints described by the polyhedra $V(s), s \in S$.

We will weaken the simplicity requirement and introduce the notion of *piecewise homogeneous hidden Markov chain*.

DEFINITION 5. A process $\mathbf{z} = (\mathbf{z}_n)$ over $S$ is a (piecewise homogeneous) hidden Markov chain if there are a finite set $S'$ and a (piecewise homogeneous) Markov chain $\mathbf{w} = (\mathbf{w}_n)$ over $S \times S'$ such that $\mathbf{z}$ is the projection of $\mathbf{w}$ over $S$.

Thus, a hidden Markov chain is the projection of a Markov chain with values in a product space. Correspondingly, a piecewise homogeneous hidden Markov chain is the projection of a piecewise homogeneous Markov chain.

We are now in position to describe the problem considered in this section. Given a sequence $x = (x_0, \ldots, x_N)$ in $S$ with finite length $N+1$, does there exist a stochastic process $\mathbf{z}$ that (i) is both a $V$-process and a piecewise homogeneous hidden Markov chain, and (ii) approximates $x$ in the sense that both assertions (B1) and (B2) in Theorem 1 hold?

Without further qualifications, the answer is negative. Indeed, if all $V$-processes are transient, assertion (B1) cannot hold. On the other hand, if the sequence $x$ is not *typical*, in the sense that the observed transition function $p^x$ is far from $V$ [i.e., $p^x(s,\cdot)$ is far from $V(s)$ in the Euclidean norm for some $s \in S$], assertion (B2) cannot hold. The following two examples illustrate these points. In both examples $V(s)$ is a singleton for each $s \in S$, hence there is a unique $V$-process which is a Markov chain.

EXAMPLE (A nonirreducible Markov chain). Let $S = \{a, b, c\}$. Define $V$ so that both states $b$ and $c$ are absorbing, while state $a$ leads with equal probability to states $b$ and $c$. Starting from state $a$, one of the two sequences $(a, b, b, b, \ldots, b)$ and $(a, c, c, c, \ldots, c)$ results. But if the given sequence is, for example, $x = (a, b, b, b, \ldots, b)$, the unique $V$-process does not satisfy (B1) when starting from state $a$.

We shall therefore restrict our study to sets $V$ such that there exists an irreducible homogeneous $V$-Markov chain.

EXAMPLE (A nontypical sequence). Let $S = \{a, b\}$ and define $V$ so that both states lead with equal probability to $a$ and $b$. If the given sequence is $x = (a, a, \ldots, a)$ the unique $V$-process does not satisfy (B2).

We shall therefore limit ourselves to sequences that are typical w.r.t. $V$, in the following sense:

DEFINITION 6. Let $N \in \mathbf{N}$ and $\delta, \varepsilon > 0$ be given. A sequence $x = (x_0, \ldots, x_N)$ is $(N, \delta, \varepsilon)$-*typical* if there exists $v = (v(s, \cdot))_s \in V$ such that $|1 - \frac{v(s,t)}{p^x(s,t)}| < \varepsilon$ for every $s, t \in S$ such that either $N_s^x p^x(s,t) \geq N^\delta$ or $N_s^x v(s,t) \geq N^\delta$.

In Definition 6, $N_s^x$ is the observed number of stages spent in $s$ along $x$, and $p^x$ is the observed transition function; see Section 2.1.

In words, a sequence $x$ is typical if there is a transition function $v$ such that (i) $v(s) \in V(s)$ for every $s$: $v$ is an admissible transition, and (ii) $v(s,t)$ is close to $p^x(s,t)$ whenever the transition from $s$ to $t$ occurs frequently, either along $x$ or under $v$. The latter statement is not completely accurate, as we do not use the invariant distribution of $v$, but rather the observed occupancy measure.

The set of $(N, \delta, \varepsilon)$-typical sequences is denoted by $\mathcal{T}(N, \delta, \varepsilon)$. Note that the notion of a typical sequence is relative to the family $V$ of polyhedra.

Our first theorem states that if $V$ contains an irreducible transition function, and if $x$ is typical w.r.t. $V$, then one can approximate $x$ in a proper sense by a piecewise homogeneous hidden Markov chain.

THEOREM 4. *Assume that $V$ contains an irreducible transition function $b$, and set $B = \max_{s,t \in S} \mathbf{E}_{s,b}[T_t]$. Let $\psi, \eta \in (0,1)$ be given. There exist $\delta, \varepsilon > 0$ and $N_1 \in \mathbf{N}$ such that the following holds. For every $N \geq N_1$ and every sequence $x \in \mathcal{T}(N, \delta, \varepsilon)$, there exists a $V$-piecewise homogeneous hidden Markov chain $\mathbf{z}$, with at most $|S|$ pieces, such that the following hold:*

(G1) $\mathbf{P}_\mathbf{z}(|\frac{\nu_N^\mathbf{z}(s)}{\nu^x(s)} - 1| \geq \eta) \leq \frac{1}{N^\delta}$ *for each $s \in S$ such that $\nu^x(s) \geq \frac{1}{N^\delta}$.*

(G2) *Let $\mathbf{N}_0 = |\{n < N : \|p_n^\mathbf{z}(\mathbf{z}_n, \cdot) - p^x(\mathbf{z}_n, \cdot)\|_\infty > \eta\}|$. Then $\mathbf{E}_\mathbf{z}[\mathbf{N}_0] \leq N^\psi B$.*

Our second theorem states that if $x$ is generated by a $V$-process, then with high probability it is typical.

THEOREM 5. *Let $\delta, \varepsilon > 0$ and $\xi \in (0, \delta/4)$ be given. There exists $N_2 \in \mathbf{N}$ such that, for every $N \geq N_2$ and every $V$-process $\mathbf{z}$, one has*

$$\mathbf{P}_\mathbf{z}(\mathcal{T}(N, \delta, \varepsilon)) \geq 1 - \frac{1}{N^\xi}.$$

These two theorems can be combined as follows. Let us postulate that we get to observe *some* realization of *some* $V$-process. Then, with high probability, we will be able to find a simple $V$-process that typically yields the observed realization.

This section is organized as follows. In Section 3.2 we prove Theorem 5, and we then turn to the proof of Theorem 4 in Section 3.3.

3.2. *Typical sequences*: *proof of Theorem 5.* We here prove Theorem 5. The proof uses the following large deviation estimate for Bernoulli variables. Let $(X_n)_{n \in \mathbf{N}}$ be an infinite sequence of i.i.d. Bernoulli r.v.s with parameter $p$, and set $\overline{X}_n = \sum_{i=1}^n X_i/n$, for each $n \in \mathbf{N}$. By Alon, Spencer and Erdös [(2000), Corollary A.14],

$$\mathbf{P}(|\overline{X}_n - p| > \varepsilon p) \leq 2\exp(-c_\varepsilon pn),$$

where $c_\varepsilon = \min\{\varepsilon^2, -\varepsilon + (1+\varepsilon)\ln(1+\varepsilon)\}$ is independent of $n$ and $p$. Hence, for $k \in \mathbf{N}$,

$$(26) \quad \mathbf{P}\left(\sup_{pn \geq k} |\overline{X}_n - p| > \varepsilon p\right) \leq 2 \sum_{n=\lceil k/p \rceil}^\infty \exp(-c_\varepsilon pn) \leq \frac{2\exp(-c_\varepsilon k)}{1 - \exp(-c_\varepsilon p)}.$$

Observe that for every $\varepsilon$ sufficiently small, $\varepsilon^2/3 < c_\varepsilon \leq \varepsilon^2/2$.

Let $\delta, \varepsilon \in (0,1)$ and $\xi \in (0, \delta/4)$ be given. Choose $\xi' \in (\xi, \delta/4)$ and set $\varepsilon' = \frac{\varepsilon}{(1+\varepsilon)\max_{s \in S}|V^*(s)| + \varepsilon}$. Let $N_2 \in \mathbf{N}$ be large enough so that the following inequalities are satisfied for each $N \geq N_2$: (T1) $\frac{2\exp(-c_{\varepsilon'} N^{\delta/4})}{1 - \exp(-c_{\varepsilon'} N^{\delta/4 - 1})} \leq 1/N^{\xi'}$, (T2) $N^{\xi' - \xi} \geq 3|S|^2 \sum_{s \in S}|V^*(s)|$, and (T3) $N^{\delta/2} \geq \frac{1 - \varepsilon'}{\varepsilon'}$.

Let $N \geq N_2$, and let $\mathbf{z} = (\mathbf{z}_n)_n$ be a $V$-process. We start by introducing a convenient decomposition for $\mathbf{z}$. Recall that for $s \in S$, $V^*(s)$ is the finite set of the extreme points of $V(s)$. For each $n \geq 0$, the conditional distribution of $\mathbf{z}_{n+1}$, given $\mathbf{z}_0, \ldots, \mathbf{z}_n$, belongs to $V(\mathbf{z}_n)$, hence can be written as a convex combination $\sum_{v \in V^*(\mathbf{z}_n)} \mathbf{b}_n(v) v$, where the weights $\mathbf{b}_n(v)$ are random. We then divide the choice of $\mathbf{z}_{n+1}$ into two steps. In the first step a point $\mathbf{v}_n \in V^*(\mathbf{z}_n)$ is drawn according to the weights $\mathbf{b}_n$. Next $\mathbf{z}_{n+1}$ is drawn according to $\mathbf{v}_n$. In other words, we simply view the given $V$-process as a process $\mathbf{z} = ((\mathbf{z}_n, \mathbf{v}_n))_n$, where $\mathbf{v}_n \in V^*(\mathbf{z}_n)$ and the conditional law of $\mathbf{z}_{n+1}$, given past values, is $\mathbf{v}_n$, for each $n \geq 0$.

Following the notation used in Section 2, we let $\mathbf{N}^\mathbf{z}_{s,v,t} = |\{n < N, (\mathbf{z}_n, \mathbf{v}_n, \mathbf{z}_{n+1}) = (s,v,t)\}|$ denote the number of one-step transitions from $s$ to $t$ which occurred through the extreme point $v$. We also use a number of derived notions. For instance, $\mathbf{N}^\mathbf{z}_{s,\cdot,t} = \sum_{v \in V^*(s)} \mathbf{N}^\mathbf{z}_{s,v,t}$ is the number of one-step transitions from $s$ to $t$, $\mathbf{N}^\mathbf{z}_{s,v,\cdot} = \sum_{t \in S} \mathbf{N}^\mathbf{z}_{s,v,t}$ is the number of visits to $s$ that are followed by the choice of the extreme point $v$, $\mathbf{N}^\mathbf{z}_s = \sum_{v \in V^*(s), t \in S} \mathbf{N}^\mathbf{z}_{s,v,t}$ is the number of visits to $s$, and $\mathbf{p}^\mathbf{z}((s,v),t) = \mathbf{N}^\mathbf{z}_{s,v,t}/\mathbf{N}^\mathbf{z}_{s,v,\cdot}$ is the proportion of transitions to $t$, out of $(s,v)$. It is defined only when $\mathbf{N}^\mathbf{z}_{s,v,\cdot} > 0$. Note that the empirical transition function is given by

$$\mathbf{p}^\mathbf{z}(s,t) = \frac{\mathbf{N}^\mathbf{z}_{s,\cdot,t}}{\mathbf{N}^\mathbf{z}_s} = \frac{\sum_{v \in V^*(s)} \mathbf{N}^\mathbf{z}_{s,v,\cdot} \mathbf{p}^\mathbf{z}((s,v),t)}{\mathbf{N}^\mathbf{z}_s}.$$

Finally, we set $\mathbf{u}^\mathbf{z}(s,\cdot) = \frac{\sum_{v \in V^*(s)} \mathbf{N}^\mathbf{z}_{s,v,\cdot} v}{\mathbf{N}^\mathbf{z}_s}$. It is the average point of $V(s)$ used at $s$. Since $V(s)$ is convex, $\mathbf{u}^\mathbf{z}(s,\cdot) \in V(s)$ for each $s \in S$. We will show that with high $\mathbf{P}_\mathbf{z}$-probability, $\mathbf{u}^\mathbf{z}$ is close to $\mathbf{p}^\mathbf{z}$ in the sense of Definition 6.

We now move to the core of the argument. The following lemma asserts that if the transition $(s,v) \mapsto t$ occurs frequently, then with high probability, the observed probability $\mathbf{p}^{\mathbf{z}}((s,v),t)$ of moving from $(s,v)$ to $t$ is close to the true one, $v(t)$.

LEMMA 6. *Let $s,t \in S$ and $v \in V^*(s)$ be given. Then*

$$\mathbf{P_z}\left(\mathbf{N}^{\mathbf{z}}_{s,v,\cdot}\max\{v(t), \mathbf{p}^{\mathbf{z}}((s,v),t)\} \geq N^{\delta/2} \Rightarrow \left|\frac{\mathbf{p}^{\mathbf{z}}((s,v),t)}{v(t)} - 1\right| \leq \varepsilon'\right) \tag{27}$$

$$\geq 1 - \frac{3}{N^{\xi'}}.$$

PROOF. Note first that $\mathbf{N}^{\mathbf{z}}_{s,v,\cdot}v(t) < N^{\delta/4}$ if $v(t) < N^{\delta/4-1}$. Assume now that $v(t) \geq N^{\delta/4-1}$. Let $(X_n)_{n \leq N}$ be a sequence of i.i.d. Bernoulli r.v.'s with parameter $v(t)$. By (26) and (T1),

$$\mathbf{P_z}(\mathbf{N}^{\mathbf{z}}_{s,v,\cdot}v(t) \geq N^{\delta/4} \text{ and } |\mathbf{p}^{\mathbf{z}}((s,v),t) - v(t)| > \varepsilon'v(t)) \leq \frac{1}{N^{\xi'}}. \tag{28}$$

Moreover, one has

$$\mathbf{P_z}(\mathbf{N}^{\mathbf{z}}_{s,v,\cdot}v(t) < N^{\delta/4} \text{ and } \mathbf{N}^{\mathbf{z}}_{s,v,\cdot}\mathbf{p}^{\mathbf{z}}((s,v),t) \geq N^{\delta/2}) \leq 2/N^{\delta/4}. \tag{29}$$

Indeed, let $(X_i)$ be a sequence of i.i.d. Bernoulli r.v.s with parameter $v(t)$, and set $n = \lfloor N^{\delta/4}/v(t) \rfloor$. By Markov's inequality, the left-hand side in (29) is at most

$$\mathbf{P_z}\left(\max_{k\,:\,kv(t)<N^{\delta/4}}\{X_1 + \cdots + X_k\} \geq N^{\delta/2}\right) \leq \mathbf{P_z}(X_1 + \cdots + X_n \geq N^{\delta/2})$$

$$\leq nv(t)/N^{\delta/2} \leq 2/N^{\delta/4}.$$

To conclude, equation (27) follows from (28), (29) and the choice of $\xi'$.

Let $\mathcal{T}^*$ be the set of all realizations $y = (x_0, v_0, x_1, v_1, \ldots, x_N)$ for which the implication in (27) holds, for every $s,t \in S$ and every $v \in V^*(s)$. By Lemma 6 and (T2), $\mathbf{P_z}(\mathcal{T}^*) \geq 1 - \frac{3|S|^2 \sum_{s \in S}|V^*(s)|}{N^{\xi'}} \geq 1 - \frac{1}{N^{\xi}}$. To conclude the proof of Theorem 5, it therefore suffices to show that every sequence in $\mathcal{T}^*$ is $(N,\delta,\varepsilon)$-typical.

Let $y$ be a sequence in $\mathcal{T}^*$. Following earlier use, we denote by $N^y_s$ the value of $\mathbf{N}^{\mathbf{z}}_s$ at $y$, and we use a similar convention for other random variables. We shall verify that $y \in \mathcal{T}(N,\delta,\varepsilon)$.

Let $s,t \in S$ be such that $N^y_s p^y(s,t) \geq N^\delta$. We will prove that $|1 - p^y(s,t)/u^y(s,t)| \leq \varepsilon$. The argument is also valid if $N^y_s u^y(s,t) \geq N^\delta$. It is enough to prove that

$$N^y_{s,v,\cdot}|p^y((s,v),t) - v(s,t)| \leq \frac{\varepsilon'}{1-\varepsilon'}N^y_s p^y(s,t) \qquad \text{for every } v \in V^*(s). \tag{30}$$

Indeed, by summing (30) over all $v \in V^*(s)$, it follows that

$$|p^y(s,t) - u^y(s,t)| \leq \sum_{v \in V^*(s)} \frac{N^y_{s,v,\cdot}}{N^y_s}|p^y((s,v),t) - v(t)| \leq \frac{\varepsilon'}{1-\varepsilon'}p^y(s,t)|V^*(s)|,$$

which implies, by the choice of $\varepsilon'$, that $|p^y(s,t) - u^y(s,t)| \leq \varepsilon u^y(s,t)$, as desired.

We let $v \in V^*(s)$ be given and proceed to the proof of (30). If $N^y_{s,v,\cdot} \max\{p^y((s,v),t), v(t)\} \geq N^{\delta/2}$, then, since the implication in (27) holds for $y$, one has $|v(t) - p^y((s,v),t)| \leq \varepsilon' v(t)$. Multiplying both sides by $N^y_{s,v,\cdot}$, we get

$$N^y_{s,v,\cdot}|p^y((s,v),t) - v(t)| \leq \varepsilon' N^y_{s,v,\cdot} v(t)$$
$$\leq \frac{\varepsilon'}{1-\varepsilon'} N^y_{s,v,\cdot} p^y((s,v),t) \leq \frac{\varepsilon'}{1-\varepsilon'} N^y_s p^y(s,t),$$

where the last inequality holds since $N^y_s p^y(s,t) = \sum_{v \in V^*(s)} N^y_{s,v,\cdot} p^y((s,v),t)$, and (30) holds. If, on the other hand, $N^y_{s,v,\cdot} \max\{p^y((s,v),t), v(t)\} < N^{\delta/2}$, then, since $N^y_s p^y(s,t) \geq N^\delta$,

$$N^y_{s,v,\cdot}|p^y((s,v),t) - v(t)| \leq N^{\delta/2} \leq N^y_s p^y(s,t)/N^{\delta/2},$$

and (30) holds by (T3). □

3.3. *Proof of Theorem* 4. We first provide a heuristic overview of the proof. It will be helpful to contrast it with the proof given in Section 2. In the basic setup, the given sequence $x$, or equivalently, the given array $(N^x_{s,t})_{s,t \in S}$ of one-step transitions, was first extended to a periodic and exhaustive sequence. Next, the structure theorem was used to find a certain partition into atoms. The approximating process simply visited each atom in turn for a number of stages equal to the observed one. The transition function of the process was such that each atom was a recurrent set. It was obtained by watching the observed transition function on the different atoms. Moving from one atom to another was done in a single step. These last features allowed for a simple analysis.

At a broad level, the analysis of the general problem is similar. We again start by extending the given sequence $x$ to a periodic and exhaustive sequence $x^*$, and by applying the structure theorem to obtain a partition $S_1, \ldots, S_k$ of $S$ (see Section 3.3.2). As in Section 2, the approximating process will focus on each atom in turn.

However, here we are constrained to use $V$-processes, hence, the former choice of a transition function may not be feasible. Instead, we introduce the transition function $v \in V$ that is closest to $p^{x^*}$ (we omit details in this sketch). The approximating process will essentially evolve according to $v$. To be more precise, consider a specific phase $k$. Since $S_k$ need not be recurrent for $v$, the process may occasionally exit from $S_k$. We will then let it evolve according to $b$, so as to re-enter $S_k$ in a few stages (recall that $b \in V$ is a fixed irreducible transition function). In a first approximation, the behavior of the approximating process during phase $k$ can thus be described by the transition function that coincides with $v$ on $S_k$, and with $b$ on $\overline{S}_k$.

It turns out that it is convenient to amend this definition as follows. Once the process exits from $S_k$, a (fictitious) entry state $t$ in $S_k$ is drawn, according to the distribution of the entry state under a specific Markov chain (again, we omit details). The process will evolve according to $b$ until $t$ is reached. It then switches back to $v$. Note that this no longer describes a Markov chain, since the transition function on $S_k$ may be either $b$ or $v$, depending upon the circumstances. This feature is best dealt with by adding a component in the state space, which keeps track of the current status of the process. This component takes values in $S' = S \cup \{\circ\}$, where $\circ$ is an additional symbol. The $k$th piece of the approximating process is defined as the $S$-marginal of a Markov chain over $S \times S'$, whose transition function $\pi_k$ is defined as follows. Whenever the $S'$-component is set to $\circ$, the $S$-component evolves according to $v$. As long as the $S$-component remains in $S_k$, the $S'$-component remains equal to $\circ$. When the $S$-component exits $S_k$, then an element $t$ of $S_k$ is selected with a given probability, and the $S'$-component of the Markov chain is set to $t$. This $t$ is the target of the $S$-component, which evolves according to $b$ as long as $t$ is not reached. Once $t$ is reached, the $S'$-component is set to $\circ$. For the purpose of the transition from phase $k-1$ to phase $k$, the exact definition of the transition function will be slightly different.

The new aspects raise additional difficulties.

First, note that the set $\mathbf{N}_0$ that appears in the statement (G2) roughly coincides with the set of stages in which the process moves according to $b$. In order to prove that the cardinality of this set is small compared to $N$, one needs to prove that the expected time to reach $S_k$ under $b$ is small compared to the expected time to leave $S_k$ under $v$. The expected time to leave $S_k$ under $p^{x^*}$ can be derived from the sequence $x^*$. We will thus have to compare the expected exit times from $S_k$, computed under $v$ and $p^{x^*}$. To do that, we will use a result on the comparison of exit times from a given set under close Markov chains.

Second, in order to prove (G1), we need to compare the empirical frequency for which the $S$-component is $s \in S_k$ with the frequency of $s$ along $x$. To this end, we prove, as in Section 2, that the transition function that is defined by $v$ on $S_k$ and by $b$ on $\overline{S}_k$ is mixing. We then use a result relative to close Markov chains in order to compare the invariant distribution of the latter transition function to the invariant distribution of $p^{x^*}$.

We now describe the organization of the proof. We define the approximating process **z** in Sections 3.3.2 and 3.3.3. We then state in Section 3.3.4 two propositions that readily imply Theorem 4. These two propositions are statements about the transition function $\pi_k$. The subsequent sections are devoted to the proofs of these propositions. Sections 3.3.5 and 3.3.6 contain the statement and the application to our framework of results on perturbed Markov chains. These results are used in the last three sections, which conclude the proof.

3.3.1. *Fixing parameters.* Let $\psi, \eta \in (0, 1)$ be given. We here list a number of conditions on $\varepsilon, \delta$ and $N_1$ under which the conclusion of Theorem 4 holds. We stress that we do not strive for optimal conditions.

Fix $0 < \varepsilon < \eta/56L < \eta$, with $L = \sum_{n=1}^{|S|} \binom{|S|}{n} n^{|S|}$.

Choose $\beta \in (0, \frac{1}{2}(\frac{A}{L})^{|S|} \times \frac{\varepsilon(1-\varepsilon)}{BL^2 \times |S|^4})$, where $A = 1/2$. Set $\alpha = \frac{1}{2\beta|S|L^2}$ and $\alpha' = \frac{\alpha/2 - |S|}{2|S|}$. Note that $\beta < 1/20|S|^2 L^2$, so that $\alpha' \geq 2$. Choose $\psi' \in (0, \psi)$, $\xi \in (0, \psi'/(|S|+1))$, $\delta' \in (0, \xi/2)$. Finally, choose $\delta \in (0, \min\{\delta', (1-\psi)/2\})$. Set $a = N^\xi$.

Choose $N_1 \in \mathbf{N}$ sufficiently large such that (C1)–(C8) and (A1)–(A7) hold, for every $N \geq N_1$: (C1) $N^{\delta'} \geq N^\delta + 1$, (C2) $N^\delta \geq \frac{3}{\varepsilon(1-\varepsilon)}$, (C3) $2 + 8L|S|(N + |S|)^{\psi'} \times (1 + |S|^2/N^\delta) \leq N^\psi/|S|$, (C4) $N^{\xi-\delta} \geq 1/\eta$, (C5) $N^{1-\delta-\psi'} \geq 8BL|S|$, (C6) $N^{1-2\delta-\psi} \geq 42(B+1)|S|/\varepsilon^2$, (C7) $\varepsilon N^{1+\xi-\delta} \geq 2(N+|S|)$, (C8) $N^\delta \geq |S|^2(1+55\varepsilon L)/\varepsilon$, (A1) $\beta(N^\xi - 1) \geq (N+|S|)^{\delta'}$, (A2) $N^\xi - 1 \geq \frac{1}{2\beta|S|}$, (A3) $L(N^\xi + 1)^{|S|} \leq N_*^{\psi'}$, (A4) $\frac{17 \times 8}{\varepsilon^2 N^{1-\delta}}(\frac{N}{N^\xi - 1} + B + 1) \leq \frac{N^{-\delta}}{L|S|^2}$, (A5) $B(1+3\varepsilon)\frac{(N^\xi+1)^{|S|}}{N^{1-\delta}} \leq \frac{1}{2N^\delta} \leq \varepsilon$, (A6) $N^\psi/|S| \geq 1 + 2(1+3\varepsilon)N^\delta(N^\xi+1)^{|S|}$, (A7) $N^\xi \geq 18|S|$.

We will prove that the conclusion of Theorem 4 holds for every $N \geq N_1$ and every typical $S$-valued sequence $x$ of length $N + 1$.

3.3.2. *The periodicized sequence.* Let $x$ be a $(N, \delta, \varepsilon)$-typical sequence. Let $x^* = (x_0^*, \ldots, x_{N_*}^*)$ be the periodic and exhaustive sequence that is obtained by extending $x$ as we did in Section 2.4.3.

Since $x$ is typical, we can choose once and for all, for every $s \in S$, an element $v(s, \cdot) \in V(s)$ such that, for every $t \in S$,

$$(31) \quad N_s^x \max\{p^x(s,t), v(s,t)\} \geq N^\delta \implies \left|1 - \frac{v(s,t)}{p^x(s,t)}\right| \leq \varepsilon.$$

As the next lemma asserts, since $x$ is typical, so is $x^*$.

LEMMA 7. $x^*$ is $(N_*, \delta', 3\varepsilon)$-typical.

PROOF. By (C1),

$$N_s^{x^*} \max\{p^{x^*}(s,t), v(s,t)\} \geq N_*^{\delta'} \implies N_s^x \max\{p^x(s,t), v(s,t)\} \geq N^\delta.$$

In that case, by (C2), $|1 - \frac{N_{s,t}^x}{N_{s,t}^{x^*}}| \leq \frac{\varepsilon}{3}$ and $|1 - \frac{N_s^x}{N_s^{x^*}}| \leq \frac{\varepsilon}{3}$. This implies that $|1 - \frac{p^{x^*}(s,t)}{p^x(s,t)}| = |1 - \frac{N_{s,t}^{x^*}/N_s^{x^*}}{N_{s,t}^x/N_s^x}| \leq \varepsilon$. Together with (31), we deduce that

$$(32) \quad N_s^{x^*} \max\{p^{x^*}(s,t), v(s,t)\} \geq N_*^{\delta'} \implies \left|1 - \frac{v(s,t)}{p^{x^*}(s,t)}\right| \leq 3\varepsilon.$$

Thus, the extended sequence $x^*$ is $(N_*, \delta', 3\varepsilon)$-typical. □

The following lemma asserts that for every state $s$ that is frequently visited, the observed transition out of $s$, $p^x(s,\cdot)$ and the transitions out of $s$ under the $V$-Markov chain $v$, are close.

LEMMA 8. *Let $s \in S$ be given. If $N_s^x \geq N^\xi$, then $\|v(s,\cdot) - p^x(s,\cdot)t\|_\infty \leq \eta$.*

PROOF. Let $t \in S$ be given. If $\max\{v(s,t), p^x(s,t)\} \leq \eta$, then clearly $|v(s,t) - p^x(s,t)| \leq \eta$. Otherwise, by (C4), $N_s^x \max\{v(s,t), p^x(s,t)\} \geq \eta N^\xi \geq N^\delta$. Therefore, by (31) and the choice of $\varepsilon$, $|v(s,t) - p^x(s,t)| \leq \varepsilon p^x(s,t) < \eta$. □

3.3.3. *The approximating process.* We here construct the approximating hidden Markov chain **z**. Its properties will be established in later sections.

We apply Theorem 2 to the sequence $x^*$ and $a = N^\xi$, and obtain a partition $\mathcal{C} = (S_1, S_2, \ldots, S_K)$ of $S$. Let $K_0 = \{k : N_{S_k}^{x^*} \geq N^{1-\delta}\}$ be the frequently visited atoms. For convenience, we assume that $K_0$ consists of the first $|K_0|$ atoms of the partition $\mathcal{C}$, so that $K_0 = \{1, \ldots, |K_0|\}$. We assume also that $S_{|K_0|}$ is the most frequently visited atom, so that, in particular, $N_{S_{K_0}}^{x^*} \geq N/|S|$.

For $k \in K_0$, we define a transition function $\pi_k$ over $\Omega = S \times (S \cup \{\circ\})$ as follows:

1. *At state $(s, \circ)$, $s \in S_k$.* A state $s' \in S$ is first drawn according to $v(s, \cdot)$. If $s' \in S_k$, the chain moves to $(s', \circ)$; if $s' \notin S_k$, a state $t \in S_k$ is drawn according to $\mathbf{P}_{s', p^{x^*}}(T_{S_k} = T_t)$ and the chain moves to $(s', t)$.
2. *At state $(s,t)$, $s \neq t$ and $t \in S_k$.* A state $s' \in S$ is first drawn according to $b(s, \cdot)$. If $s' = t$, the chain moves to $(s', \circ)$. Otherwise, the chain moves to $(s', t)$.
3. *At state $(s,t)$, $s \notin S_k$ and $t \in \overline{S}_k \cup \{\circ\}$.* A pair $(s', t') \in \Omega$ is drawn with probability $b(s, s') \times \mathbf{P}_{s, p^{x^*}}(T_{S_k} = T_{t'})$. If $s' = t'$, the chain moves to $(s', \circ)$. Otherwise, the chain moves to $(s', t')$.

Other states are visited with probability 0. Note that the $S$-marginal of $\pi_k((s,t), \cdot)$ is either $v(s, \cdot)$ or $b(s, \cdot)$ and, in particular, belongs to $V(s)$.

Plainly, a chain with transition function $\pi_k$ always moves in a single step to a state in $(S \times S_k) \cup (S_k \times \{\circ\})$. In particular, the third item in the definition of $\pi_k$ may possibly be relevant only at the initial stage. Note that the $S$-coordinate behaves under $\pi_k$ as described in the overview. Starting from $S_k \times \{\circ\}$, it moves according to $v$, unless it exits $S_k$. In that case, a target state in $S_k$ is drawn according to the distribution of the entry state computed with $p^{x^*}$. Then the $S$-coordinate moves according to $b$ until it reaches the target state. At this point, the target flag is removed, and the $S$-coordinate resumes moving according to $v$.

The approximating hidden Markov chain has $|K_0|$ pieces. The length of piece $k$ is $N_{S_k}^{x^*}$, except that of piece $|K_0|$: its length is $N - \sum_{k<|K_0|} N_{S_k}^{x^*}$, which is between $N_{S_{|K_0|}}^{x^*}$ and $N_{S_{|K_0|}}^{x^*} + |S|N^{1-\delta}$. In piece $k$ the process follows the transition function $\pi_k$. We denote by $m_k$ the length of piece $k$.

Formally, we let $\mathbf{w}$ be a piecewise homogeneous Markov chain over $\Omega$ whose transition function coincides with $\pi_k$ at stages $\sum_{j<k} m_j \leq n < \sum_{j\leq k} m_j$. We define $\mathbf{z}$ to be the first component of $\mathbf{w}$, so that it is a piecewise hidden Markov chain. Thus, for every stage $n$ in piece $k$, the conditional law of $\mathbf{w}_{n+1}$ is $\pi_k(\mathbf{w}_n, \cdot)$. The initial state of $\mathbf{w}$ is irrelevant. We will prove that the process $\mathbf{z}$ satisfies both (G1) and (G2).

For the convenience of the proof, the definition of the boundaries of the $k$th piece slightly differs from the one in Section 2.4.2.

3.3.4. *Two propositions.* We here state two propositions relative to $\pi_k$, without proof. We next show why Theorem 4 follows from these propositions. As a consequence, the proof of Theorem 4 reduces to statements about Markov chains.

In Proposition 2 below, $\nu_{m_k}(s, \circ)$ is the empirical frequency of visits to the state $(s, \circ)$ in stages 0 through $m_k - 1$. Recall that $\nu^{x^*}(s|S_k) = \nu^{x^*}(s)/\nu^{x^*}(S_k)$.

PROPOSITION 2.  *Let $k \in K_0$ be given. For every $\omega \in \Omega$ and every $s \in S_k$,*

$$\mathbf{P}_{\omega, \pi_k}\left( \left| 1 - \frac{\nu_{m_k}(s, \circ)}{\nu^{x^*}(s|S_k)} \right| > 55\varepsilon L \right) \leq \frac{1}{2N^\delta}.$$

In effect, Proposition 2 contains two statements. By summation over $s \in S_k$, it implies that, with high probability, the empirical frequency of $S_k \times \{\circ\}$ is close to one. It also says that the empirical frequency of $(s, \circ)$ is close to the observed frequency of $s$ along the sequence $x^*$, when conditioned on $S_k$.

In Proposition 3 below, $\mathbf{N}_{\overline{S_k \times \{\circ\}}} = |\{n < m_k \text{ s.t. } \mathbf{w}_n \notin (S_k \times \{\circ\})\}|$ is the number of stages in which the $S$-coordinate of the state evolves according to $b$ rather than according to $v$. Recall that $B = \max_{s,t \in S} \mathbf{E}_{s,b}[T_t]$ is a bound on the expected time under $b$ to reach any given state.

PROPOSITION 3.  *Let $k \in K_0$ and $\omega \in \Omega$ be given. One has*

$$\mathbf{E}_{\omega, \pi_k}[\mathbf{N}_{\overline{S_k \times \{\circ\}}}] \leq \frac{1}{|S|} B N^\psi.$$

In effect, Proposition 3 asserts that the total number of stages the process spends outside $S_k \times \{\circ\}$ is small. As a consequence, the empirical frequency of $S_k \times \{\circ\}$ is close to one. This statement, however, differs from Proposition 2, since here the result is phrased in expected terms.

We now show why the conclusions of Theorem 4 follow from Propositions 2 and 3.

We begin with (G2) and let $m_{k-1} \leq n < m_k$. For every $s \in S_k$, one has, by Theorem 2(P2), $N_s^x \geq R_{\{s\}}^x > N^\xi$. Hence, by Lemma 8, $|v(s,\cdot) - p^x(s,\cdot)| \leq \eta$. By definition of $\pi_k$, the $S$-marginal of $p_n^\mathbf{w}(\mathbf{w}_n, \cdot)$ is equal to $v(\mathbf{z}_n, \cdot)$ whenever $\mathbf{w}_n \in S_k \times \{\circ\}$. In particular, $\mathbf{w}_n \in S_k \times \{\circ\}$ implies $n \notin \mathbf{N}_0$. Denote by $\mathbf{N}_{0,k}$ the number of visits to states outside $S_k \times \{\circ\}$ during the $k$th piece. By construction, $\mathbf{E}_\mathbf{z}[\mathbf{N}_{0,k}] \leq \sup_{\omega \in \Omega} \mathbf{E}_{\omega, \pi_k}[\mathbf{N}_{\overline{S_k \times \{\circ\}}}]$, so that by Proposition 3,

$$(33) \qquad \mathbf{E}_\mathbf{z}[\mathbf{N}_0] \leq \sum_{k \in K_0} \mathbf{E}_\mathbf{z}[\mathbf{N}_{0,k}] \leq BN^\psi,$$

and (G2) follows.

We next check that (G1) holds. Fix $s \in S$, such that $\nu^x(s) \geq 1/N^\delta$. By construction $s \in S_k$ for some $k \in K_0$. We introduce the frequency $\tilde{\nu}_{m_k}^\mathbf{w}(s, \circ) = \frac{1}{m_k}|\{m_k \leq n < m_{k+1} : \mathbf{w}_n = (s, \circ)\}|$ of visits to $(s, \circ)$ in piece $k$. Note that the difference $N\nu_N^\mathbf{z}(s) - m_k \tilde{\nu}_{m_k}^\mathbf{w}(s, \circ)$ is the sum of two terms: (i) the total number of visits to $s$ in phases other than $k$, and (ii) the total number of visits to $\{s\} \times S$ during phase $k$. As a consequence, $N\nu_N^\mathbf{z}(s) - m_k \tilde{\nu}_{m_k}^\mathbf{w}(s, \circ) \leq \sum_{k \in K_0} \mathbf{N}_{0,k}$. By (33), Markov's inequality and (C6), one has

$$(34) \qquad \mathbf{P}_\mathbf{w}\left(N\nu_N^\mathbf{z}(s) - m_k \tilde{\nu}_{m_k}^\mathbf{w}(s, \circ) > \varepsilon \frac{N}{N^\delta}\right) \leq \frac{BN^{\delta+\psi-1}}{\varepsilon} \leq \frac{1}{2N^\delta}.$$

Note that the conditional distribution of $\tilde{\nu}_{m_k}^\mathbf{w}(s, \circ)$, given $\mathbf{w}_0, \ldots, \mathbf{w}_{m_k}$, coincides with the distribution of $\nu_{m_k}(s, \circ)$ under a Markov chain starting from $\mathbf{w}_{m_k}$ and with transition $\pi_k$. Hence, by Proposition 2,

$$(35) \qquad \mathbf{P}_\mathbf{w}\left(\left|1 - \frac{\tilde{\nu}_{m_k}^\mathbf{w}(s, \circ)}{\nu^{x^*}(s|S_k)}\right| > 55\varepsilon L\right) \leq \frac{1}{2N^\delta}.$$

By (34) and (35), the probability that both inequalities $N\nu_N^\mathbf{z}(s) - m_k \tilde{\nu}_{m_k}^\mathbf{w}(s, \circ) \leq \varepsilon N^{1-\delta}$ and $|1 - \frac{\tilde{\nu}_{m_k}^\mathbf{w}(s, \circ)}{\nu^{x^*}(s|S_k)}| \leq 55\varepsilon L$ hold is at least $1 - \frac{1}{N^\delta}$. On this event, by (C8),

$$|\nu_N^\mathbf{z}(s) - \nu^{x^*}(s)|$$
$$\leq \left|\nu_N^\mathbf{z}(s) - \frac{m_k}{N}\tilde{\nu}_{m_k}^\mathbf{w}(s, \circ)\right| + \frac{1}{N}\tilde{\nu}_{m_k}^\mathbf{w}(s, \circ)|m_k - N_{S_k}^{x^*}|$$
$$+ \frac{N_{S_k}^{x^*}}{N}|\tilde{\nu}_{m_k}^\mathbf{w}(s, \circ) - \nu^{x^*}(s|S_k)|$$
$$\leq \varepsilon N^{-\delta} + \frac{1}{N}|S|^2 N^{-\delta} N_{S_k}^{x^*} \tilde{\nu}_{m_k}^\mathbf{w}(s, \circ) + \frac{N_{S_k}^{x^*}}{N} \times 55\varepsilon L \nu^{x^*}(s|S_k)$$

$$\leq \varepsilon \frac{N+|S|}{N}\nu^{x^*}(s) + \frac{N_{S_k}^{x^*}}{N}|S|^2 N^{-\delta}(1+55\varepsilon L)\nu^{x^*}(s|S_k) + 55\varepsilon L\nu^{x^*}(s)$$
$$\leq 56\varepsilon L\nu^{x^*}(s).$$

By the choice of $\varepsilon$, this proves (G1).

3.3.5. *Perturbation of Markov chains*: *reminder.* We here introduce a result on perturbations of Markov chains due to Solan and Vieille (2003). This result provides an estimate of the sensitivity of the stationary distribution and other statistical quantities with respect to perturbations of the transition function. This result will be applied to our setup in the next section.

Given $C \subseteq S$ with $|C| \geq 2$, and an irreducible transition function $p^1$ over $S$ with invariant measure $\mu^1$, set

$$\zeta_{p^1}^C = \min_{\varnothing \subset D \subset C} \sum_{s \in D} \mu^1(s)p^1(s,\overline{D}).$$

This is a variation of the conductance of a Markov chain that was originally defined by Jerrum and Sinclair (1989) and was used in the study of the rate of convergence to the invariant measure [see, e.g., Lovász and Kannan (1999) and Lovász and Simonovits (1990)].

The notion of closeness we use is the following one:

DEFINITION 7. Let $p^1$ be an irreducible transition function on $S$ with invariant measure $\mu^1$, let $C \subseteq S$ with $|C| \geq 2$, and let $\beta, \varepsilon > 0$. A transition function $p^2$ is $(\beta,\varepsilon)$-*close to* $p^1$ *on* $C$ if (i) $p^2(s,\cdot) = p^1(s,\cdot)$ for every $s \notin C$, and (ii) $|1 - \frac{p^2(s,t)}{p^1(s,t)}| \leq \varepsilon$ for every $s,t \in C$ such that $\mu^1(s)\max\{p^1(s,t),p^2(s,t)\} \geq \beta\zeta_{p^1}^C$.

This definition is not symmetric since it involves the invariant distribution of $p^1$, and not that of $p^2$. It requires that the relative probabilities of moving from $s$ to $t$ under $p^1$ and under $p^2$ are close, but only for those one-step transitions $s \to t$ that, on average, occur relatively frequently.

The next result summarizes Theorems 4 and 6 in Solan and Vieille (2003). It asserts that if $p^1$ and $p^2$ are two transition functions that are close in the sense of Definition 7, then their invariant measures are close, as well as other statistical quantities of interest, such as the average length of a visit to a set and the exit time from a given set. Recall that $L = \sum_{n=1}^{|S|-1}\binom{|S|}{n}n^{|S|}$ and that $T_C$ is the first hitting time of the set $C$. The quantity $K_p(C)$ has been introduced in Section 2.3.2.

PROPOSITION 4. *Let $\varepsilon \in (0, 1/2^{|S|})$, $A > 0$, and $\beta \in (0, \frac{1}{2}(\frac{A}{L})^{|S|} \times \frac{\varepsilon(1-\varepsilon)}{L \times |S|^4})$ be given. Let $p^1$ be an irreducible transition function defined over $S$ with invariant measure $\mu^1$. Assume that $|C| \geq 2$ and that $\mathbf{P}_{s,p^1}(T_t^+ < T_{\overline{C}}) \geq A$ for every $s, t \in C$. Let $p^2$ be $(\beta, \varepsilon)$-close to $p^1$ on $C$. Then:*

(a) *All states of $C$ belong to the same recurrent set for $p^2$. Let $\mu^2$ be the invariant measure of $p^2$ on that recurrent set.*

(b) *For every $s \in C$ and every $D \subset C$, one has*

$$\left|1 - \frac{\mu^2(s|C)}{\mu^1(s|C)}\right| \leq 18\varepsilon L, \tag{36}$$

$$L^{-1} \leq \frac{\mathbf{E}_{s,p^2}[T_{\overline{D}}]}{\mathbf{E}_{s,p^1}[T_{\overline{D}}]} \leq L \quad \text{and} \quad L^{-1} \leq \frac{K_{p^2}(D)}{K_{p^1}(D)} \leq L. \tag{37}$$

(c) *Let $\chi \in (0, \beta\zeta_{p^1}^C]$ be any number such that, for every $s, t \in C$,*

$$\mu^1(s) \max\{p^1(s,t), p^2(s,t)\} \geq \chi \implies \left|1 - \frac{p^2(s,t)}{p^1(s,t)}\right| \leq \varepsilon. \tag{38}$$

*Then either*

$$L^{-1} K_{p^1}(C) \leq K_{p^2}(C) \leq L K_{p^1}(C), \quad \text{or} \tag{39}$$

$$K_{p^1}(C) \geq \frac{1}{2|S|} \times \frac{\mu^1(C)}{\chi} \quad \text{and} \quad K_{p^2}(C) \geq \frac{1}{L} \times \frac{1}{2|S|} \times \frac{\mu^1(C)}{\chi}. \tag{40}$$

3.3.6. *Perturbation of Markov chains: application.* We here introduce the auxiliary transition function $q_k$ on $S$ defined by

$$q_k = \begin{cases} v & \text{on } S_k, \\ p^{x^*} & \text{on } \overline{S_k}. \end{cases}$$

In Lemmas 9 and 10 below, we first check that the conditions of Proposition 4 are fulfilled by $q_k$ and $p^{x^*}$, as soon as $S_k$ is not a singleton. Next, relying on Proposition 4, we provide estimates of the mixing measures $\lambda_{q_k}(S_k)$ and $\rho_{q_k}(S_k)$ (see Proposition 5 below). These estimates will later be used to relate the properties of $q_k$ to those of the transition function $\pi_k$ over $\Omega$.

LEMMA 9. *If $|S_k| \geq 2$, the transition function $q_k$ is $(\beta, 3\varepsilon)$-close to $p^{x^*}$ on $S_k$.*

PROOF. By Definition 7 and (32), it suffices to prove that $\beta\zeta_{p^{x^*}}^{S_k} \geq \frac{N_*^{\delta'}}{N_*}$. For each $C \subset S_k$, one has by (P2) and since $a = N^\xi$,

$$\sum_{s \in C} \nu^{x^*}(s) p^{x^*}(s, \overline{C}) = \frac{R_C^{x^*} - \mathbb{1}_{x_N \in C}}{N_*} \geq \frac{N^\xi - 1}{N_*}. \tag{41}$$

By taking the minimum over $C \subset S_k$, this yields $\zeta_{p^{x^*}}^{S_k} \geq \frac{N^\xi - 1}{N_*}$. The result follows by (A1). □

Recall that we set $A = 1/2$.

LEMMA 10. *If $|S_k| \geq 2$, one has $\mathbf{P}_{s,p^{x^*}}(T_t^+ < T_{\overline{S}_k}) \geq A$ for every $s, t \in S_k$.*

PROOF. Suppose first that $s \neq t$, so that $\mathbf{P}_{s,p^{x^*}}(T_t^+ < T_{\overline{S}_k}) = \mathbf{P}_{s,p^{x^*}}(T_t < T_{\overline{S}_k})$. By Lemma 2,

$$\mathbf{P}_{s,p^{x^*}}(T_t < T_{\overline{S}_k}) \geq 1 - 2|S_k| \frac{\rho_{p^{x^*}}(S_k)}{\lambda_{p^{x^*}}(S_k) - (|S_k| - 1)\rho_{p^{x^*}}(S_k)}.$$

By Lemma 5, $\rho_{p^{s^*}}(S_k) \leq \frac{2}{N^\xi} \lambda_{p^{x^*}}(S_k)$, so that by (A7),

(42) $\quad \mathbf{P}_{s,p^{x^*}}(T_t < T_{\overline{S}_k}) \geq 1 - 2|S_k| \frac{2/N^\xi}{1 - 2(|S_k| - 1)/N^\xi} \geq \frac{3}{4} \geq A.$

Suppose now that $t = s$. Since $|S_k| \geq 2$, $p^{x^*}(s, \overline{S}_k) \leq \frac{R_{S_k}^{x^*}}{R_s^{x^*}} < \frac{1}{a} = \frac{1}{N^\xi}$, so that by (42) and (A7),

$$\mathbf{P}_{s,p^{x^*}}(T_s^+ > T_{\overline{S}_k}) = p^{x^*}(s, \overline{S}_k) + \sum_{t \in S_k \setminus \{s\}} p^{x^*}(s,t) \mathbf{P}_{t,p^{x^*}}(T_s > T_{\overline{S}_k})$$

$$\leq \frac{1}{N^\xi} + \frac{1}{4} \leq \frac{1}{2} = 1 - A. \qquad \square$$

By Lemmas 9 and 10, we can apply Proposition 4 to $p^{x^*}$ and $q_k$ with $A = 1/2$. Recall that $\alpha = 1/(2\beta |S| L^2)$.

PROPOSITION 5. *Assume that $|S_k| \geq 2$. Then $\lambda_{q_k}(S_k) \geq \alpha \rho_{q_k}(S_k)$ and $\lambda_{q_k}(S_k) \geq \frac{1}{2L|S|} \frac{N_{S_k}^{x^*}}{N_*^{\psi'}}.$*

PROOF. We first provide a lower bound on $K_{p^{x^*}}(S_k)$. By (16), one has

(43) $\quad \frac{K_{p^{x^*}}(S_k)}{\nu^{x^*}(S_k)} = \frac{1}{\sum_{s \in S_k} \nu^{x^*}(s) p^{x^*}(s, \overline{S}_k)} = \frac{N_*}{N_{S_k, \overline{S}_k}^{x^*}} \geq \frac{N_*}{R_{S_k}^{x^*}}.$

We now prove the first assertion. For $C \subset S_k$, by (43) and (P2),

$$K_{p^{x^*}}(S_k) \geq \frac{N_* \nu^{x^*}(S_k)}{R_{S_k}^{x^*}} \geq N^\xi \frac{N_*}{R_C^{x^*}} \nu^{x^*}(S_k) \geq (N^\xi - 1) \frac{\nu^{x^*}(S_k)}{\sum_{s \in C} \nu^{x^*}(s) p^{x^*}(s, \overline{C})},$$

where the last inequality follows since $\frac{N^\xi}{N_{C,\bar{C}}+1} \geq \frac{N^\xi - 1}{N_{C,\bar{C}}}$. By optimizing over $C$ and using (A2), this yields

$$(44) \qquad K_{p^{x^*}}(S_k) \geq \frac{1}{2\beta|S|} \times \frac{\nu^{x^*}(S_k)}{\zeta_{p^{x^*}}^{S_k}}.$$

Using (39) and (40) with $\chi = \beta \zeta_{p^{x^*}}^{S_k}$, (44) yields

$$(45) \qquad \lambda_{q_k}(S_k) \geq K_{q_k}(S_k) \geq \frac{1}{L} \times \frac{1}{2\beta|S|} \times \frac{\nu^{x^*}(S_k)}{\zeta_{p^{x^*}}^{S_k}}.$$

Fix $C \subset S_k$. By (45) and (37),

$$\lambda_{q_k}(S_k) \geq \frac{1}{2\beta|S|L} \times \frac{\nu^{x^*}(C)}{\sum_{s \in C} \nu^{x^*}(s) p^{x^*}(s, \overline{C})} = \frac{1}{2\beta|S|L} \times K_{p^{x^*}}(C)$$

$$\geq \frac{1}{L} \times \frac{1}{2\beta|S|L} \times K_{q_k}(C) \geq \frac{1}{L^2} \times \frac{1}{2\beta|S|} \times \min_{s \in C} \mathbf{E}_{s, q_k}[T_{\overline{C}}].$$

The first assertion follows by taking the maximum over $C$.

We now prove the second assertion. By (32) and since $\delta' < \psi'$, (38) holds with $\chi = N_*^{\psi'-1}$. We distinguish two cases. If $K_{p^{x^*}}(S_k) \geq \frac{1}{2|S|} \times \frac{N_{S_k}^{x^*}}{N_* \times \chi}$, then by (40),

$$\lambda_{q_k}(S_k) \geq K_{q_k}(S_k) \geq \frac{1}{L} \times \frac{1}{2|S|} \times \frac{N_{S_k}^{x^*}}{N_* \times \chi} = \frac{1}{L} \times \frac{1}{2|S|} \times \frac{N_{S_k}^{x^*}}{N_*^{\psi'}},$$

as desired. If, on the other hand, $K_{p^{x^*}}(S_k) < \frac{1}{2|S|} \frac{N_{S_k}^{x^*}}{N_* \times \chi}$, then by (39) and (43),

$$\lambda_{q_k}(S_k) \geq K_{q_k}(S_k) \geq \frac{1}{L} K_{p^{x^*}}(S_k)$$

$$\geq \frac{1}{L} \times \frac{N_* \nu^{x^*}(S_k)}{R_{S_k}^{x^*}} \geq \frac{1}{L} \times \frac{N_{S_k}^{x^*}}{R_{S_k}^{x^*}} \geq \frac{N_{S_k}^{x^*}}{L(N^\xi + 1)^{|S|}},$$

which by (A3) gives the result. □

3.3.7. *Proof of Proposition* 3 *when* $|S_k| \geq 2$. By Proposition 4, there is a recurrent set for $q_k$ that contains $S_k$. Therefore, there is a recurrent set $\Omega_k \subset \Omega$ for $\pi_k$ that contains $S_k \times \{\circ\}$. We denote by $\mu_{\pi_k}$ the invariant measure of $\pi_k$ on $\Omega_k$. We first prove that $\mu_{\pi_k}$ assigns a significant weight to $S_k \times \{\circ\}$.

LEMMA 11. *If* $|S_k| \geq 2$, *then* $\mu_{\pi_k}(S_k \times \{\circ\}) \geq 1 - \frac{2B}{\lambda_{q_k}(S_k)} \geq \frac{1}{2}$.

PROOF. By Proposition 5 and (C5), $1 - \frac{2B}{\lambda_{q_k}(S_k)} \geq \frac{1}{2}$.

We now prove that $\mu_{\pi_k}(S_k \times \{\circ\}) \geq 1 - \frac{2B}{\lambda_{q_k}(S_k)}$. Plainly, $\mathbf{E}_{(s,t),\pi_k}[T_{S_k \times \{\circ\}}] = \mathbf{E}_{s,b}[T_t] \leq B$ for every $t \in S_k$ and every $s \in S \setminus \{t\}$.

On the other hand, by Lemma 1(ii), Proposition 5 and by the choice of $\beta$,

$$\min_{s \in S_k} \mathbf{E}_{(s,\circ),\pi_k}[T_{\overline{S_k \times \{\circ\}}}] = \min_{s \in S_k} \mathbf{E}_{s,q_k}[T_{\overline{S_k}}]$$
$$\geq \lambda_{q_k}(S_k) - (|S_k| - 1)\rho_{q_k}(S_k) \geq \tfrac{1}{2}\lambda_{q_k}(S_k).$$

By (18), one gets

$$\frac{\mu_{\pi_k}(\overline{S_k \times \{\circ\}})}{\mu_{\pi_k}(S_k \times \{\circ\})} \leq \frac{2B}{\lambda_{q_k}(S_k)},$$

hence $\mu_{\pi_k}(\overline{S_k \times \{\circ\}}) \leq \frac{2B}{\lambda_{q_k}(S_k)}$.

We now proceed to the proof of Proposition 3 when $|S_k| \geq 2$. Observe first that

$$\min_{\omega \in S_k \times \{\circ\}} \mathbf{E}_{\omega,\pi_k}[\mathbf{N}_{\overline{S_k \times \{\circ\}}}] = \min_{\omega \in \Omega_k} \mathbf{E}_{\omega,\pi_k}[\mathbf{N}_{\overline{S_k \times \{\circ\}}}].$$

Since $\mu_{\pi_k}$ is the invariant measure of $\pi_k$ over $\Omega_k$, one has $\mathbf{E}_{\mu_{\pi_k},\pi_k}[\mathbf{N}_{\overline{S_k \times \{\circ\}}}] = m_k \mu_{\pi_k}(\Omega_k \setminus (S_k \times \{\circ\}))$. By Lemma 11 this yields

$$(46) \qquad \min_{\omega \in S_k \times \{\circ\}} \mathbf{E}_{\omega,\pi_k}[\mathbf{N}_{\overline{S_k \times \{\circ\}}}] \leq 2B \times \frac{m_k}{\lambda_{q_k}(S_k)}.$$

Let $\gamma = \max_{\omega \in S_k \times \{\circ\}} \mathbf{E}_{\omega,\pi_k}[\mathbf{N}_{\overline{S_k \times \{\circ\}}}]$, and let $\omega_1 \in S_k \times \{\circ\}$ achieve the maximum. Since the $S$-marginal of $\pi_k$ coincides with $b$ outside $S_k \times \{\circ\}$, one has, for $\omega \in S_k \times \{\circ\}$,

$$\gamma = \mathbf{E}_{\omega_1,\pi_k}[\mathbf{N}_{\overline{S_k \times \{\circ\}}}] \leq \mathbf{E}_{\omega,\pi_k}[\mathbf{N}_{\overline{S_k \times \{\circ\}}}] + \mathbf{P}_{\omega_1,\pi_k}(T_{\overline{S_k \times \{\circ\}}} < T_\omega)(B + \gamma).$$

By Lemma 2 and Proposition 5, $\mathbf{P}_{\omega_1,\pi_k}(T_{\overline{S_k \times \{\circ\}}} < T_\omega) \leq \frac{2|S|}{\alpha/2 - |S|} = 1/\alpha'$. Since $\alpha' \geq 2$, one gets, by letting $\omega$ vary,

$$(47) \qquad \begin{aligned} \gamma &\leq \frac{\alpha'}{\alpha' - 1} \min_{\omega \in S_k \times \{\circ\}} \mathbf{E}_{\omega,\pi_k}[\mathbf{N}_{\overline{S_k \times \{\circ\}}}] + \frac{B}{\alpha' - 1} \\ &\leq 2 \min_{\omega \in S_k \times \{\circ\}} \mathbf{E}_{\omega,\pi_k}[\mathbf{N}_{\overline{S_k \times \{\circ\}}}] + B. \end{aligned}$$

Finally, for each $\omega' \in \Omega$, by (47), (46), and Proposition 5,

$$\mathbf{E}_{\omega',\pi_k}[\mathbf{N}_{\overline{S_k \times \{\circ\}}}] \leq \mathbf{E}_{\omega',\pi_k}[T_{S_k \times \{\circ\}}] + \max_{\omega \in S_k \times \{\circ\}} \mathbf{E}_{\omega,\pi_k}[\mathbf{N}_{\overline{S_k \times \{\circ\}}}]$$
$$\leq B + B + 2 \min_{\omega \in S_k \times \{\circ\}} \mathbf{E}_{\omega,\pi_k}[\mathbf{N}_{\overline{S_k \times \{\circ\}}}]$$

$$\leq 2B + 4B \times \frac{m_k}{\lambda_{q_k}(S_k)}$$

$$\leq 2B + 8BL|S|N_*^{\psi'} \times \frac{m_k}{N_{S_k}^{x^*}}.$$

Since $m_k/N_{S_k}^{x^*}$ is either 1 (if $k < |K_0|$) or at most $1 + |S|^2/N^\delta$ (if $k = |K_0|$), the desired result follows by (C3). $\square$

3.3.8. *Proof of Proposition 2 when $|S_k| \geq 2$.* To prove Proposition 2 when $|S_k| \geq 2$, we first prove that $\pi_k$ is mixing (see Lemma 13). We can then apply Theorem 3 as we did in the proof of Theorem 1. We are therefore able to compare the empirical frequency $\nu_{m_k}$ to $\mu_{\pi_k}$. Since $p^{x^*}$ and $q_k$ are close, this enables us to compare the invariant measure of $\pi_k$, $\mu_{\pi_k}$, to the invariant measure of $p^{x^*}$, $\nu^{x^*}$.

LEMMA 12. *If $|S_k| \geq 2$, then, for every $\omega \in \Omega_k$ and every $s \in S_k$, one has*

$$\mathbf{E}_{\omega,\pi_k}[T^+_{(s,\circ)}] \leq \frac{(|S_k| - 1)\rho_{q_k}(S_k) + 2B}{\min_{t \in S_k} \mathbf{P}_{t,v}(T_s < T_{\overline{S}_k})} + 1.$$

PROOF. The proof is a simple adaptation of the proof of Lemma 3. We repeat it, with few modifications. Let $\omega \in \Omega_k$ and $s \in S_k$ be given. Note that

(48) $$\mathbf{E}_{\omega,\pi_k}[T^+_{(s,\circ)}] \leq 1 + \max_{\omega' \in \Omega_k} \mathbf{E}_{\omega',\pi_k}[T_{(s,\circ)}].$$

Set $\gamma = \max_{t \in S_k} \mathbf{E}_{(t,\circ),\pi_k}[T_{(s,\circ)}]$. Let $t \in S_k \setminus \{s\}$ achieve the maximum in the definition of $\gamma$. By Lemma 1,

$$\gamma = \mathbf{E}_{(t,\circ),\pi_k}[T_{(s,\circ)}] \leq \mathbf{E}_{t,q_k}[T_{\overline{S}_k \cup \{s\}}] + \mathbf{P}_{t,q_k}(T_{\overline{S}_k} < T_s)(\gamma + B)$$
$$\leq (|S_k| - 1)\rho_{q_k}(S_k) + B + \gamma \times \max_{u \in S_k} \mathbf{P}_{u,v}(T_{\overline{S}_k} < T_s).$$

Therefore

(49) $$\gamma \leq \frac{(|S_k| - 1)\rho_{q_k}(S_k) + B}{\min_{u \in S_k} \mathbf{P}_{u,v}(T_s < T_{\overline{S}_k})}.$$

For $\omega \in \Omega_k \setminus (S_k \times \{\circ\})$,

(50) $$\mathbf{E}_{\omega,p_k}[T_{(s,\circ)}] \leq B + \gamma.$$

The result follows from (48)–(50). $\square$

The lemma below is a mixing-type result. It is very similar to Lemma 1.

LEMMA 13. *If $|S_k| \geq 2$, then, for every $\omega \in \Omega_k$ and every $s \in S_k$,*
$$\mathbf{E}_{\omega,\pi_k}[T^+_{(s,\circ)}] \leq 2|S|L\frac{N_*}{N^\xi - 1} + 4B + 1.$$

PROOF. We repeat the proof of Proposition 1 with minor adjustments. By Lemma 12,
$$\mathbf{E}_{\omega,\pi_k}[T^+_{(s,\circ)}] \leq \frac{(|S_k| - 1)\rho_{q_k}(S_k) + 2B}{\min_{t \in S_k} \mathbf{P}_{t,v}(T_s < T_{\overline{S}_k})} + 1.$$

By Lemma 2, the denominator is at least $1 - 2|S_k|\frac{\rho_{q_k}(S_k)}{\lambda_{q_k}(S_k) - (|S_k|-1)\rho_{q_k}(S_k)} \geq \frac{1}{2}$, where the inequality follows by Proposition 5 and the choice of $\beta$. Therefore
$$\mathbf{E}_{\omega,\pi_k}[T^+_{(s,\circ)}] \leq 2|S_k|\rho_{q_k}(S_k) + 4B + 1.$$

By (37), Lemma 5 and Theorem 2(P2),
$$\rho_{q_k}(S_k) \leq L\rho_{p^{x^*}}(S_k) \leq L \max_{D \subset S_k} \frac{N_D^{x^*}}{R_D^{x^*} - 1} \leq L\frac{N_*}{a - 1}.$$

The result follows. □

Define $\mu^\circ_{\pi_k}(s) = \mu_{\pi_k}((s,\circ))/\mu_{\pi_k}(S_k \times \{\circ\})$. It is the invariant measure of $\pi_k$ conditioned on $S_k \times \{\circ\}$.

The following compares the empirical number of visits to $(s, \circ)$ to the invariant measure.

PROPOSITION 6. *If $|S_k| \geq 2$, then*
$$\mathbf{P}_{\omega,\pi_k}\left(|\nu_{m_k}(s,\circ) - \mu^\circ_{\pi_k}(s)| > \frac{2\varepsilon}{1-\varepsilon}\mu^\circ_{\pi_k}(s) + \frac{1}{m_k}\right) \leq \frac{1}{2N^\delta}.$$

PROOF. By Remarks 2 and 3, Lemma 13 and (A4), for every $\omega \in \Omega_k$,
$$\mathbf{P}_{\omega,\pi_k}\left(|\nu_{m_k}(s,\circ) - \mu_{\pi_k}((s,\circ))| > \varepsilon\mu_{\pi_k}((s,\circ)) + \frac{1}{m_k}\right)$$
$$(51) \qquad \leq \frac{17}{\varepsilon^2 m_k}\left(2|S|L\frac{N_*}{N^\xi - 1} + 4B + 1\right) \leq \frac{1}{2N^\delta}.$$

Since $\mu_{\pi_k}((s,\circ)) \geq 1/2$, by Lemma 11, Proposition 5 and by the choice of $\beta$,
$$|\mu_{\pi_k}((s,\circ)) - \mu^\circ_{\pi_k}(s)| \leq 2\mu_{\pi_k}(\Omega_k \setminus (S_k \times \{\circ\})) \times \mu_{\pi_k}(s,\circ)$$
$$\leq \frac{4B}{\lambda_{q_k}(S_k)} \times \mu_{\pi_k}(s,\circ) \leq \varepsilon\mu_{\pi_k}(s,\circ).$$

Therefore, if $|\nu_{m_k}(s, \circ) - \mu_{\pi_k}((s, \circ))| \leq \varepsilon \mu_{\pi_k}((s, \circ)) + \frac{1}{m_k}$, then $|\nu_{m_k}(s, \circ) - \mu_k^\circ(s))| \leq \frac{2\varepsilon}{1-\varepsilon}\mu_k^\circ(s) + \frac{1}{m_k}$. The result follows by (51). □

We are now in a position to prove Proposition 2. Observe that the invariant measure of $p^{x^*}$ conditioned on $S_k$ is simply $\nu^{x^*}(\cdot|S_k)$. By Lemmas 9 and 10, and Proposition 4,

$$|\mu_{\pi_k}^\circ(s) - \nu^{x^*}(s \mid S_k)| \leq 18 \times 3\varepsilon L \nu^{x^*}(s|S_k).$$

The claim follows by Proposition 6, the choice of $\varepsilon$ and since by (P2) and (C7) $\nu^{x^*}(s|S_k) \geq \frac{R_{\{s\}}^{x^*}}{N_*} \geq \frac{N^\xi}{N_*} \geq \frac{2}{\varepsilon N^{1-\delta}} \geq \frac{2}{\varepsilon m_k}$.

3.3.9. *The singleton case.* We here assume that $S_k = \{s\}$ is a singleton. The next lemma is an analog of Lemma 11. It bounds from below the invariant distribution of $\pi_k$ on $S_k \times \{\circ\}$. Its proof is, however, significantly different.

LEMMA 14. *One has* $\mu_{\pi_k}((s, \circ)) \geq 1 - B(1 + 3\varepsilon)\frac{(N^\xi+1)^{|S|}}{N^{1-\delta}}$.

PROOF. By Theorem 2(P1),

$$p^{x^*}(s, S \setminus \{s\}) \leq \frac{R_{\{s\}}^{x^*}}{N_s^{x^*}} \leq \frac{(N^\xi + 1)^{|S|}}{N^{1-\delta}}.$$

Using (32), this yields

$$v(s, S \setminus \{s\}) \leq (1 + 3\varepsilon)\frac{(N^\xi + 1)^{|S|}}{N^{1-\delta}}.$$

We apply (18) to $p = \pi_k$, $S = \Omega_k$ and $C = \{(s, \circ)\}$, and we get

$$\frac{1 - \mu_{\pi_k}((s, \circ))}{\mu_{\pi_k}((s, \circ))} \leq Bv(s, S \setminus \{s\}) \leq B(1 + 3\varepsilon)\frac{(N^\xi + 1)^{|S|}}{N^{1-\delta}}.$$

The desired result follows. □

The rest of the proof for the singleton case follows closely the proof for $|S_k| \geq 2$. We first prove Proposition 2 in that case. By the definition of $\pi_k$, $\max_{\omega \in \Omega_k} \mathbf{E}_{\omega, \pi_k}[T_{(s,\circ)}^+] \leq B + 1$. Therefore, using Remark 2 with $\pi_k$, $\varepsilon$ and $\omega = (s, \circ)$, and by (C6),

$$\mathbf{P}_{\omega, \pi_k}\left(|\nu_{m_k}(s, \circ) - \mu_{\pi_k}((s, \circ))| > \varepsilon \mu_{\pi_k}((s, \circ)) + \frac{1}{m_k}\right) \leq \frac{17(B+1)}{m_k \varepsilon^2}$$
$$\leq \frac{1}{2|S|N^\delta}.$$

By Lemma 14 and (A5), $|\mu_{\pi_k}((s,\circ)) - 1| \leq \varepsilon$. The conclusion of Proposition 2 follows. Observe that we also deduce that $\mu_{\pi_k}((s,\circ)) \geq 1/2$.

We now prove Proposition 3. Fix $\omega \in \Omega_k$. Since $\mu_{\pi_k}((s,\circ)) \geq 1/2$,

$$\mathbf{E}_{\omega,\pi_k}[\mathbf{N}_{\overline{S_k \times \{\circ\}}}] \leq \mathbf{E}_{\omega,\pi_k}[T_{(s,\circ)}] + \mathbf{E}_{(s,\circ),\pi_k}[\mathbf{N}_{\overline{S_k \times \{\circ\}}}] \leq B + \mathbf{E}_{(s,\circ),\pi_k}[\mathbf{N}_{\overline{S_k \times \{\circ\}}}]$$
$$\leq B + 2\mathbf{E}_{\mu_{\pi_k},\pi_k}[\mathbf{N}_{\overline{S_k \times \{\circ\}}}] \leq B + 2m_k(1 - \mu_{\pi_k}((s,\circ))).$$

By Lemma 14 and (A6) this last quantity is at most

$$B + 2NB(1 + 3\varepsilon)\frac{(N^\xi + 1)^{|S|}}{N^{1-\delta}} \leq BN^\psi/|S|,$$

and Proposition 3 follows.

**Acknowledgments.** We thank Ehud Lehrer, Sean Meyn, Laurent Saloff-Coste and Ron Shamir for their suggestions and references. We also thank an anonymous referee whose suggestions were most helpful and substantially improved the presentation.

## REFERENCES

- ALDOUS, D. J. and FILL, J. A. (2002). Reversible Markov chains and random walks on graphs. Available at www.stat.berkeley.edu/users/aldous.
- ALON, N., SPENCER, J. H. and ERDÖS, P. (2000). *The Probabilistic Method*, 2nd ed. Wiley, New York. MR1885388
- BAUM, L. E. and EAGON, J. A. (1967). An inequality with applications to statistical estimation for probabilistic functions of Markov processes and to a model for ecology. *Bull. Amer. Math. Soc.* **73** 360–363. MR210217
- BAUM, L. E. and PETRIE, T. (1966). Statistical inference for probabilistic functions of finite state Markov chains. *Ann. Math. Statist.* **37** 1554–1563. MR202264
- FELLER, W. (1968). *An Introduction to Probability Theory and Its Applications* **1**, 3rd ed. Wiley, New York. MR228020
- GLYNN, P. W. and ORMONEIT, D. (2002). Hoeffding's inequality for uniformly ergodic Markov chains. *Statist. Probab. Lett.* **56** 143–146. MR1881167
- JERRUM, M. and SINCLAIR, A. (1989). Approximating the permanent. *SIAM J. Comput.* **18** 1149–1178. MR1025467
- KROGH, A., MIAN, I. S. and HAUSSLER, D. (1994). A hidden Markov model that finds genes in E. coli DNA. *Nucleic Acids Research* **22** 4768–4778.
- LOVÁSZ, L. and KANNAN, R. (1999). Faster mixing via average conductance. In *Proc. 31st Annual ACM Symposium on Theory of Computing* 282–287. ACM, New York. MR1798047
- LOVÁSZ, L. and SIMONOVITS, M. (1990). The mixing rate of Markov chains, an isoperimetric inequality, and computing the volume. In *Proc. 31st Annual Symposium on Foundations of Computer Science* **1** 346–354. IEEE, Piscataway, NJ. MR1150706
- RABINER, L. R. (1989). A tutorial on hidden Markov models and selected applications in speech recognition. *Proc. IEEE* **77** 257–286.
- SOLAN, E. and VIEILLE, N. (2003). Perturbed Markov chains. *J. Appl. Probab.* **40** 107–122. MR1953770

D. Rosenberg  
Laboratoire d'Analyse Geometrie  
et Applications Institut Galilée  
Université Paris Nord  
avenue Jean Baptiste Clément  
93430 Villetaneuse  
France  
e-mail: dinah@math.univ-paris13.fr

E. Solan  
MEDS Department  
Kellogg School of Management  
Northwestern University  
Evanston, Illinois 60208  
USA  
and  
School of Mathematical Sciences  
Tel Aviv University  
Tel Aviv 69978  
Israel  
e-mail: eilons@post.tau.ac.il  
e-mail: e-solan@kellogg.northwestern.edu

N. Vieille  
Département Finance et Economie  
HEC  
1 rue de la Libération  
78351 Jouy-en-Josas  
France  
e-mail: vieille@hec.fr